\documentclass[final]{amsart}

\usepackage{amssymb,amsmath,amsfonts}
\usepackage{graphicx,psfrag,epsfig,color}
\usepackage{setspace}
\usepackage{subfigure}
\usepackage[margin=1.5in]{geometry}

\DeclareMathOperator{\sech}{sech}
\newtheorem{thm}{Theorem}[section]
\newtheorem{cor}[thm]{Corollary}

\renewcommand\Re{\operatorname{\mathfrak{Re}}}
\renewcommand\Im{\operatorname{\mathfrak{Im}}}

\begin{document}

\title{An Efficient, Non-linear Stability Analysis for Detecting Pattern Formation in Reaction Diffusion Systems}

\author{William R. Holmes }
\address{Department of Mathematics, University of California Irvine, Irvine, CA. }
\email{wrholmes@uci.edu}
\thanks{{\tt Address:} Department of Mathematics, University of California Irvine, Irvine, CA. ({\tt wrholmes@uci.edu})}

\begin{abstract}

Reaction diffusion systems are often used to study pattern formation in biological systems.  However, most methods for understanding their behavior are challenging and can rarely be applied to complex systems common in biological applications.  I present a relatively simple and efficient, non-linear stability technique that greatly aids such analysis when rates of diffusion are substantially different.  This technique reduces a system of reaction diffusion equations to a system of ordinary differential equations tracking the evolution of a large amplitude, spatially localized perturbation of a homogeneous steady state.   Stability properties of this system, determined using standard bifurcation techniques and software, describe both linear and non-linear patterning regimes of the reaction diffusion system.  I describe the class of systems this method can be applied to and demonstrate its application.  Analysis of Schnakenberg and substrate inhibition models is performed to demonstrate the methods capabilities in simplified settings and show that even these simple models have non-linear patterning regimes not previously detected.  Analysis of a protein regulatory network related to chemotaxis shows its application in a more complex setting where other non-linear methods become intractable.  Predictions of this method are verified against results of numerical simulation, linear stability, and full PDE bifurcation analyses.

\end{abstract}

\keywords{reaction diffusion, pattern formation, nonlinear stability analysis, local perturbation analysis}

\subjclass[2000]{35K57, 92C15, 92C42}

\maketitle

\pagestyle{myheadings}
\thispagestyle{plain}
\markboth{W. R. Holmes}{The Local Perturbation Analysis}

\section{Introduction}

Reaction diffusion equations (RDEs) have provided a ubiquitous framework for studying pattern formation in chemical and biological systems \cite{Turing-52,Meinhardt-72,Lewis-93,Murray-02,Jilkine-07,Goehring-11}.  As a result of their sustained interest, numerous linear \cite{Turing-52}, weakly non-linear \cite{Rubinstein-99,Short-10,Rubinstein-12,Kaper-09}, and fully non-linear \cite{Iron-00,Kolokolnikov-05,Ward-02,Nishiura-12,Doelman-98,Wei-05,Mori-11} techniques for analyzing RDEs have been developed.  
Here I present an efficient, relatively simple addition to this toolbox, aimed at analyzing systems where diffusivities are substantially different.

A difference in rates of diffusion has been implicated as being of vital importance for patterning in numerous biological systems.  In the context of cell biology, some rates of diffusion are not just different but are vastly different, varying by factors of $100-1000$.  Numerous cell functions are controlled by regulators, such as ``GTPase's'' in cell motility \cite{Jilkine-07,Holmes12a,Dawes-07,Maree-06,Holmes12b,Mori-08}, ``ROP's'' in plant development \cite{Fu-01}, and ``Min'' proteins in bacterial division \cite{Huang-03,Huang-05}.  All of these regulatory proteins have fast and slow diffusing components since they exist in membrane bound and unbound states.  Analysis of these types of systems motivated the development of this method.

Here I consider a generic system of RDEs with a large diffusion disparity and highlight a useful method for understanding their linear and non-linear stability properties.  Consider 
\begin{align}\label{LPA-1}
\frac{\partial u}{\partial t}(x,t)&=f(u,v;p)+\mathbf{D_u} u_{xx},  \\
\frac{\partial v}{\partial t}(x,t)&=g(u,v;p)+\mathbf{D_v} v_{xx},  \nonumber
\end{align}
where $u$ and $v$ are vectors of slow and fast diffusing variables respectively, $\mathbf{D_u}, \mathbf{D_v}$ are diagonal matrices of diffusion coefficients, and $p$ is a vector of reaction parameters.  


The ``Local Perturbation Analysis" (LPA), is a non-linear stability technique applicable to systems of this type.  This method, originally devised by  Mar\'ee and Grieneisen \cite{Grieneisen-thesis}, is a bridge between linear and non-linear analysis methods having benefits of each.  Linear stability analysis (LSA) \cite{Turing-52} is straitforward and widely used, but is limited to providing linear information.  Non-linear methods, while more informative, 
are much more challenging, specific to the system being investigated, and rarely scale up to complex systems with many variables.  The LPA provides non-linear stability information beyond that of LSA, but is relatively simple to implement.

In contrast to LSA which probes stability of a homogeneous steady state (HSS) with respect a small amplitude, spatially extended perturbation, the LPA probes stability with respect to a spatially localized, large amplitude perturbation of the slow variable $u$.  As will be shown, when diffusion of $u$ (resp. $v$) is sufficiently slow (resp. fast) the perturbed region and broader domain evolve according to an approximate collection of ODE's on the timescale of reactions
\begin{subequations}\label{LP-sys}
\begin{align}
\frac{d u^g}{dt}(x,t)&=f(u^g,v^g;p), \label{LP-sys-a} \\
\frac{d v^g}{dt}(x,t)&=g(u^g,v^g;p), \label{LP-sys-b} \\
\frac{d u^l}{dt}(x,t)&=f(u^l,v^g;p)    \label{LP-sys-c} .
\end{align}
\end{subequations}
The variables $(u^g,v^g)$ represent ``global'' concentrations away from the perturbation and $u^l$ the concentration at the local perturbation.  Tracking the growth or decay of this perturbation provides stability information for \eqref{LPA-1}.  There are two primary benefits to this approach: 1) the large amplitude ``probe'' detects pattern formation in linearly stable parameter regimes, 2) the reduction to ODE's greatly simplifies its implementation.  

Applications of this method to biologically motivated reaction diffusion systems are found in \cite{Holmes12a,Holmes12b,Grieneisen-thesis}.  Rather than focus on a specific phenomena or biological system, my goal here is to explain the method itself.  I will describe the types of RDE's to which this method is applicable, its limitations, and the type of information that it can provide.  Well known examples of pattern forming systems are used to demonstrate its application and make direct comparisons between its predictions and results of classical methods (e.g. LSA, full PDE bifurcation, or numerical simulation).  In the context of a more complex chemotaxis related example, I also show this method: easily scales to larger systems with many variables, allows the user to gain a more complete overview of the parameter space structure than with other methods, and greatly aids investigation of both parametric and structural perturbations of a complex reaction network.

\section{Local Perturbation System Formulation}

I now proceed to show that the evolution of a spatially localized perturbation of a homogeneous steady state of Eqn.~~\eqref{LPA-1} evolves according to Eqs.~\eqref{LP-sys}.
Consider Eqn.~~\eqref{LPA-1} on the interval $[-1,1]$ with no flux boundary conditions, and $u, v$ in $\mathbb{R}^M$ and $\mathbb{R}^N$ respectively.  It is not necessary to assume all slow (respectively fast) variables have the same diffusivities, only that they can be divided into fast and slow diffusing classes.  For notation simplicity however, assume 
\begin{equation}
\mathbf{D_u}=\epsilon^2 \, \mathbb{I}, \qquad   \mathbf{D_v}=D \, \mathbb{I} ,
\end{equation}
where $\mathbb{I}$ is the properly sized identity matrix.  The central assumption will be that the three timescales defined by reaction kinetics, slow, and fast diffusion respectively are substantially different, i.e. $\epsilon^2 \ll 1 \ll D$.  Non-dimensionalization by domain size and the reaction timescale (so that $f,g \sim O(1)$) have been implicitly assumed.  Further assume this system has a HSS $(u^s(p),v^s(p))$ satisfying $f(u^s(p),v^s(p);p)=0=g(u^s(p),v^s(p);p)$.

Consider a highly localized perturbation of this steady state of the form 
\begin{align}\label{IC}
u(x,0)=u^s,  \qquad \qquad  v(x,0)=v^s,  \qquad \qquad  |x|> \sqrt{\epsilon},    \\
u(x,0)=u^p,  \qquad \qquad  v(x,0)=v^p,  \qquad \qquad  |x| < \sqrt{\epsilon},  \nonumber
\end{align}
where $(u^p,v^p)$ is $O(1)$ with respect to $\epsilon$ and $D$.  
Denote $R^l$ to be the local region $|x| < \sqrt{\epsilon}$ and $R^g$ the global region $|x|> \sqrt{\epsilon}$.

\subsection{Time and space scale separation}

To track the evolution of $(u,v)$ on these regions, different time and space scales must be considered.  A multiple timescale argument is applied to parse the role of reaction and diffusion effects on different timescales, and a boundary layer technique is used to separate relevant space scales.

The reaction, slow, and fast diffusion timescales inherent in this class of RDE's can be described by $t=O(1)$,  $t_u=\epsilon^2 t$, and $t_v = D t$, with the intermediate reaction timescale of most interest here.  Suppose $u=U(x,t,t_u,t_v) ,  v=V (x,t,t_u,t_v)$.  With the perturbation \eqref{IC} taken as an initial condition,  boundary layers on an $O(\epsilon)$ length scale are expected.  
Employing a stretched coordinate $\xi = (x - x_{layer}) / \epsilon$ near boundary layers, we come to two systems:
%
%
\begin{equation} \label{Outer-RDE}
\frac{ \partial{U} } { \partial t} = f + \epsilon^2 U_{xx},    \qquad   \frac{ \partial{V} } { \partial t} = g + D V_{xx}
\end{equation}
describing outer regions away from boundary layers and
\begin{equation} \label{Inner-RDE}
\frac{ \partial{U} } { \partial t} = f +  U_{\xi \xi},    \qquad   \frac{ \partial{V} } { \partial t} = g + \frac{D}{\epsilon^2}  V_{\xi \xi} ,
\end{equation}
describing dynamics in boundary layers.  I now describe the evolution of $(U,V)$ on the outer regions $R^{g,l}$ on the short and intermediate timescales.  Boundary layer effects and long timescale evolution are beyond the scope of this exposition.

\subsection{Evolution on the fast diffusion timescale}

Consider first the fast diffusion timescale and assume $U,V$ are described by first order perturbative expansions $U=U^0 + \epsilon U^1$, $V=V^0+ \epsilon V^1$.
Substituting $t_v=D t$ into Eqn.~\eqref{Outer-RDE} and collecting leading order terms
\begin{equation}\label{Fast-Time}
\frac{\partial U^0}{\partial t_v} =0 ,       \qquad        \frac{\partial V^0}{\partial t_v} = V^0_{xx} ,
\end{equation}
it is clear that in outer regions, $U^0=U^0(x,t,t_u)$ does not evolve due to either reaction or diffusion and $V^0$ simply spreads due to diffusion 
\begin{equation}\label{v0_exp}
V^0(x,t,t_u,t_v)=v^0(t,t_u) + \displaystyle \sum_{n=1}^{\infty} v^n(t,t_u) \exp (-(n \pi)^2 t_v) \cos( n \pi x) .
\end{equation}
So in each outer region $R^{l,g}$, $V^0$ evolves to a constant value exponentially quickly with $t_v$.

\begin{figure}[htb] 
  \centerline{ 
\includegraphics[scale=.8]{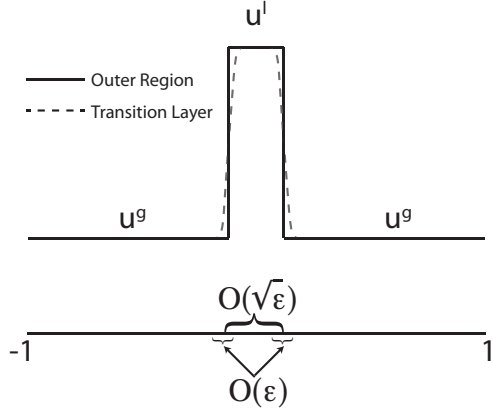}}
\caption{ To probe non-linear stability, the Local Perturbation Analysis probes the response of a homogeneous steady state of Eqn.~\eqref{LPA-1} to a localized perturbation Eqn.~~\eqref{IC}.  To leading order, values in the perturbed region and the broader domain, denoted $u^l,u^g$ respectively, will evolve independently.  The fast variable $v$, not depicted here, will be spatially uniform on the entire domain taking value $v^g$.  Boundary layer effects are present on $O(\epsilon)$ regions but to leading order do not influence evolution of the perturbation.  A collection of ordinary differential equations for $(u^g,v^g,u^l)$ describe the growth or decay of this perturbation and provide stability information for \eqref{LPA-1}.  }
\label{fig:LPA-diag}
\end{figure} 

Thus $(U^0,V^0)$ will evolve to a piecewise constant profile with possibly different values in $R^g$ and $R^l$, see Figure \ref{fig:LPA-diag}.  Denote the values on the global domain $R^g$ by $(u^g,v^g)$ and those on the local region $R^l$ as $(u^l,v^l)$.   Now consider the boundary layers between these regions on this fast timescale.  Given the spatial symmetry, consider only the left boundary layer and substitute $t_v=Dt$ into Eqn.~~\eqref{Inner-RDE}.  $O(\epsilon^{-2})$ terms indicate that $V^0_{\xi \xi}=0$.  Thus to leading order 
\begin{equation*}
V^0=a_0 (t_v) \xi + a_1 (t_v) .
\end{equation*}
Matching conditions dictate that 
\begin{equation*}
\displaystyle \lim_{\xi \to \infty} V^0 (\xi) = v^l   ,    \qquad    \displaystyle \lim_{\xi \to -\infty}  V^0 (\xi)  = v^g   .
\end{equation*}
So $a_0=0$, resulting in a shadow system \cite{Nishiura-82,Li-09} where $V^0$ is constant over the entire domain.    Its value will be denoted $v^g(t,t_u)$.  On this timescale to leading order, $U^0_{t_v}=0$ so that boundary layer effects do not influence $U^0$.

\subsection{Evolution on the intermediate reaction timescale}
Evolution of this perturbation on the reaction timescale will determine the stability of the HSS.  As $t_v$ progresses approaching this timescale, to leading order the solution in the outer regions is described by  $(u^g,v^g,u^l)$, representing the piecewise constant values on $R^{l,g}$ respectively.  The evolution of these values on this timescale are described by Eqn.~~\eqref{Outer-RDE} with $t=O(1)$.
To leading order
\begin{equation}\label{Intermediate-Time}
\frac{\partial U^0}{\partial t} =f(U^0,V^0;p) ,       \qquad        \frac{\partial V^0}{\partial t} = g(U^0,V^0;p)  .
\end{equation}
Substituting $(u^g,v^g)$ and $(u^l,v^g)$ respectively into the first of these, we obtain Eqs.~\eqref{LP-sys-a}, \eqref{LP-sys-c}.
Integrating the second of \eqref{Intermediate-Time} over the domain we see that 
\begin{equation}\label{vg-correction}
\frac{\partial v^g}{\partial t} = \frac{1}{2} \displaystyle \int_{-1}^{1} g(U^0,v^g;p) \, dx = g(u^g,v^g;p) + \sqrt{\epsilon} \, [ g(u^l,v^g;p)-g(u^g,v^g;p)] .
\end{equation}
So to leading order the values $(u^l,u^g,v^g)$ evolve according to \eqref{LP-sys}
on the intermediate timescale.  
Eqn.~\eqref{LP-sys} will be referred to as the LPA system of ODE's (or LPA-ODE's) associated with Eqn.~~\eqref{LPA-1}.  Evolution on the slow diffusion timescale requires consideration of boundary layer effects.  This is dependent on the specific system being investigated and is not discussed here.

\section{The Local Perturbation Analysis: Application and examples}\label{Important-points}

The goal of the LPA is to determine in which parameter regimes localized perturbations of this form grow or decay.  Growth corresponds to instability of the HSS and a patterning response, decay corresponds to stability.  Since the evolution of this perturbation is described by Eqn.~\eqref{LP-sys}, the location and stability of its steady states / fixed points describe stability properties of the HSS of Eqn.~\eqref{LPA-1}.  This information can be found using standard bifurcation analysis techniques for systems of ODE's.


In coming sections, I will demonstrate this method through example and the following capabilities will be emphasized.
\begin{enumerate} 
\item The LPA detects linear instabilities of Eqn.~~\eqref{LPA-1}.          

\item The LPA detects inherently non-linear patterning where a HSS is linearly stable but a sufficiently large perturbation yields a patterning response.

\item While the LPA approximation is not valid on the timescale of pattern evolution, its results can be used to make reasonable conjectures about the type of pattern (i.e. highly localized spike or a sharp interface separating distinct planer regions) that might evolve.

\item In non-linear patterning regimes, the LPA qualitatively maps the dependence of patterning response thresholds on system parameters $p$ (excluding diffusion parameters).
\end{enumerate}

The LPA is applied to two classical systems, Schnakenberg \cite{Schnakenberg-79} and Substrate Inhibition \cite{Thomas-79}, both well studied in \cite{Murray-82,Ward-02,Iron-04} for example.  Predictions of this analysis are then directly compared to results of linear stability, numerical, full PDE bifurcation, and asymptotic analyses.

\subsection{The local perturbation analysis of a Schnakenberg model}\label{Schnak-results}

The Schnakenberg system is a Turing model where $u$ is an activator and $v$ a substrate.
\begin{subequations}\label{snack-sys}
\begin{align}\label{snack-sys-u}
u_t(x,t)&=a-u+u^2v+ \epsilon^2 \triangle u=f(u,v)+ \epsilon^2 \triangle u, \\
v_t(x,t)&=b-u^2v + D \triangle v = g(u,v)+D \triangle v.   \label{snack-sys-v}
\end{align}
\end{subequations}
$u$ decays linearly, both are produced uniformly in the domain, and the nonlinearity represents an autocatalytic reaction where $u$ consumes $v$.  A LSA of Eqn.~~\eqref{snack-sys} can be found in \cite{Murray-82}.  In \cite{Ward-02,Iron-04}, asymptotic and spectral analysis showed that for $a=0$, highly localized spike solutions exist and are stable. 

It is possible to analytically perform the LPA for Eqn.~~\eqref{snack-sys}.  The resulting system of LPA-ODE's becomes 
\begin{subequations}\label{snack-LPA-sys}
\begin{align}\label{snack-LPA-ug}
u^g_t&=a-u^g+(u^{g})^{2}v^g, \\
v^g_t&=b-(u^{g})^{2} v^g,   \label{snack-LPA-vg} \\
u^l_t&=a-u^l+(u^{l})^{2}v^g .  \label{snack-LPA-ul}
\end{align}
\end{subequations}
Here, $p=(a,b)$ is the vector of system parameters and $a$ will be the bifurcation parameter of interest.
Eqs.~ \eqref{snack-LPA-ug}, \eqref{snack-LPA-vg} decouple from Eqn.~\eqref{snack-LPA-ul} and simply represent the spatially homogeneous system (i.e. with $\epsilon=0=D$).
The unique HSS $(u^s(p),v^s(p))$ of Eqn.~~\eqref{snack-sys}
\begin{equation}\label{Schnak-HSS}
 u^s=a+b,   \qquad  v^s=\frac{b}{(a+b)^2}, 
\end{equation}
is thus a solution of Eqs.~\eqref{snack-LPA-ug}, \eqref{snack-LPA-vg}.  Similarly $(u^g,v^g,u^l)=(u^s,v^s,u^s)$ is a steady state of Eqn.~~\eqref{snack-LPA-sys}.  This steady state of the LPA-ODE's represents the HSS of Eqn.~~\eqref{snack-sys} with no perturbation, i.e. both $u^{l,g}=u^s$.  

While this is the only HSS of the RDE system, the LPA system actually has two steady states with the second satisfying 
\begin{equation}\label{LPA-SS}
u^g=u^s, \qquad v^g=v^s,  \qquad  u^l=a+\frac{a^2}{b} =: u^{l1}
\end{equation}
With $b$ fixed and $a$ considered as a bifurcation parameter, the steady state branches $u^s$ and $u^{l1}$ intersect in a transcritical bifurcation at $a=b$.  Furthermore, it can be readily shown by computing the Jacobian of Eqn.~~\eqref{snack-LPA-sys} that the stability of these branches is determined solely by the sign of $f_u$ and that 
\begin{align}
\frac{\partial f}{\partial u}(u^s,v^s)&>0, & \frac{\partial f}{\partial u}(u^{l1},v^s)&<0, & a&<b, \\
\frac{\partial f}{\partial u}(u^s,v^s)&<0, & \frac{\partial f}{\partial u}(u^{l1},v^s)&>0, & a&>b. 
\end{align}

The location and stability of these steady states is depicted in Figure \ref{fig:Schnak}a.  The HSS branch $(u^g,v^g,u^l)=(u^s,v^s,u^s)$ is linearly unstable for $a<1$ (Region I).  For $a>1$ (Region II), the HSS is linearly stable, however a perturbation of $u^l$ above the $u^{l1}$ branch will grow to infinity.   From here on, the HSS branch $(u^g,v^g,u^l)=(u^s,v^s,u^s)$ will be referred to as a ``global'' steady state branch of the LPA system.  $(u^g,v^g,u^l)=(u^s,v^s,u^{l1})$ will be referred to as a ``local'' branch, since it is only a steady state for the local variable $u^l$.  

\subsubsection{Local perturbation analysis predictions}

These results lead to the following predictions.

\vspace{.1in}
\noindent \textbf{Prediction 1:} A Turing bifurcation occurs near $a=1$.  For $a<1$, the homogeneous steady state of Eqn.~~\eqref{snack-sys} is linearly unstable.  For $a>1$, it is stable but sufficiently large perturbations yield a patterning response.  
\vspace{.1in}

\begin{figure}[htb!]
\psfrag{Du}{$D_u$}
\subfigure[]{
\centering
\includegraphics[scale=.23]{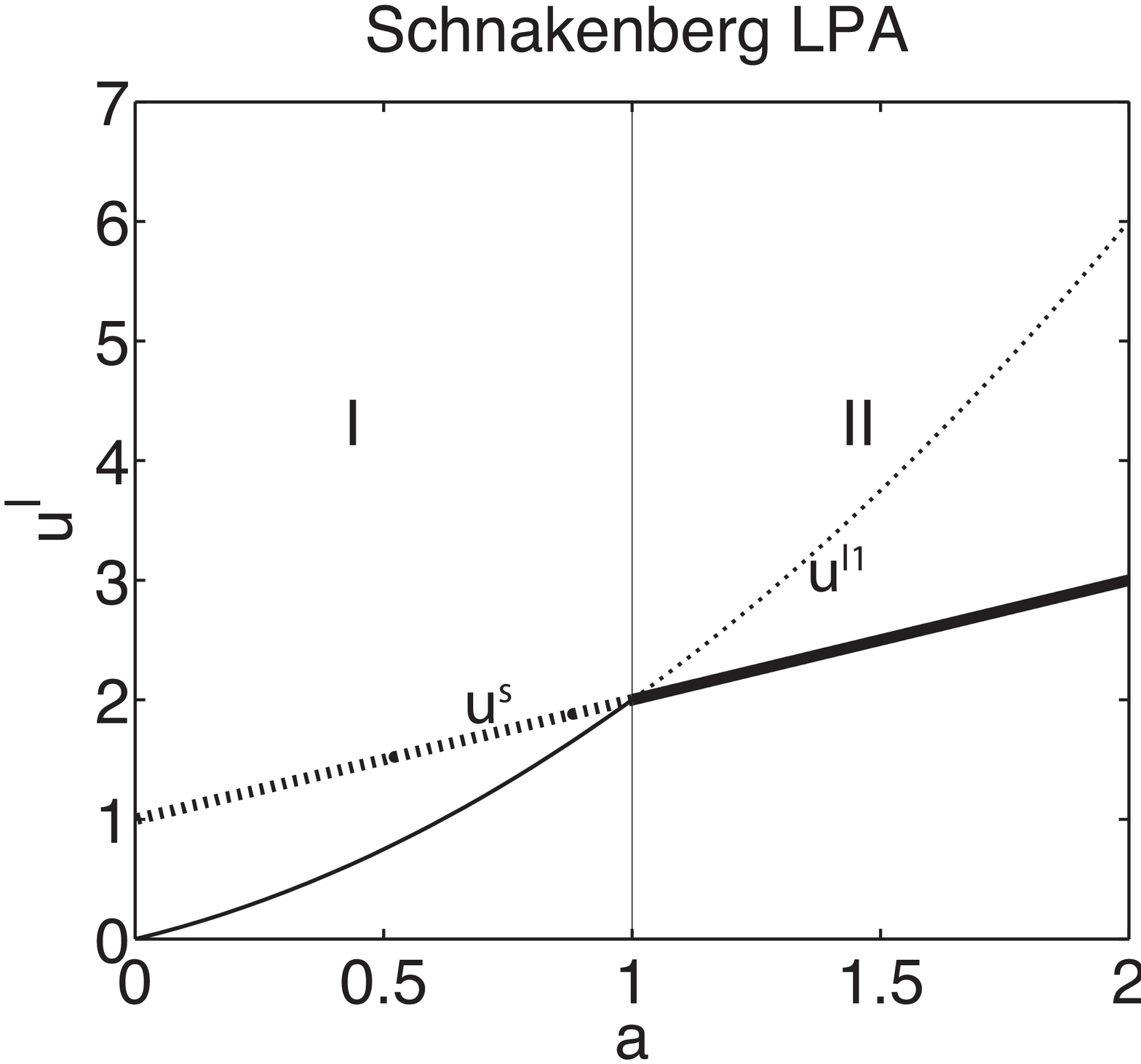}
}
\subfigure[]{
\centering
\includegraphics[scale=.23]{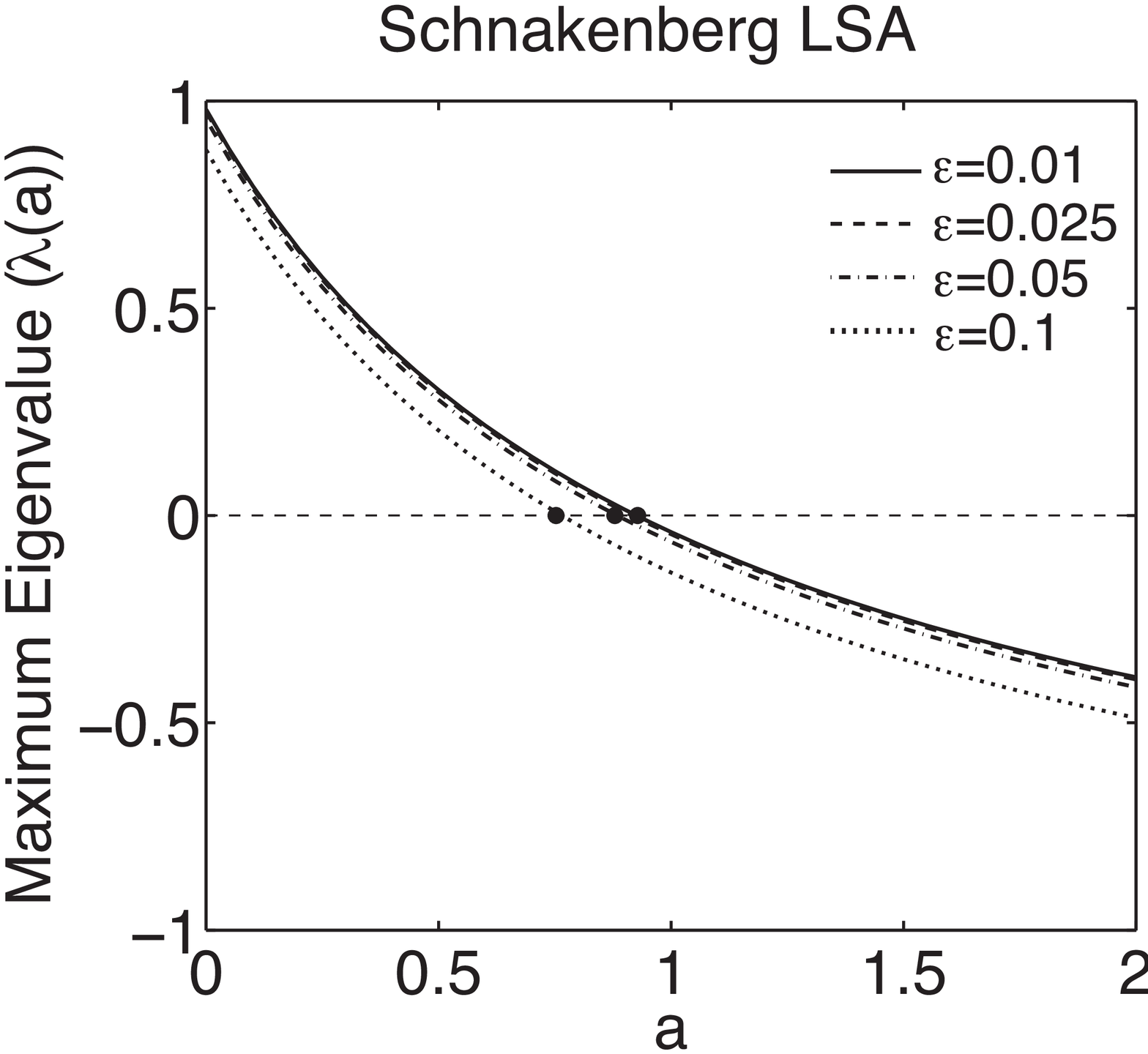}
}
\subfigure[]{
\centering
\includegraphics[scale=.23]{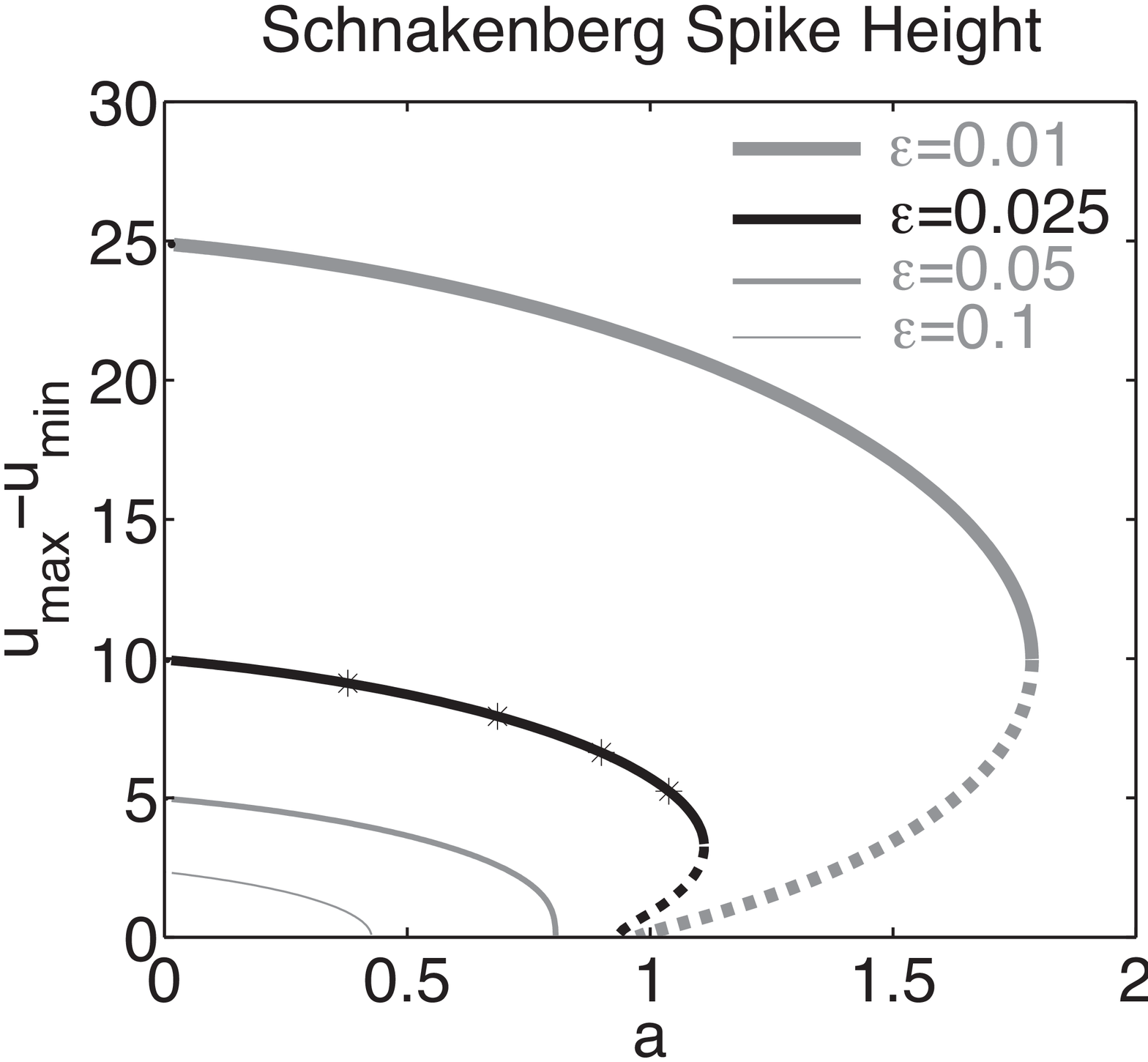}
}
\subfigure[]{
\centering
\includegraphics[scale=.23]{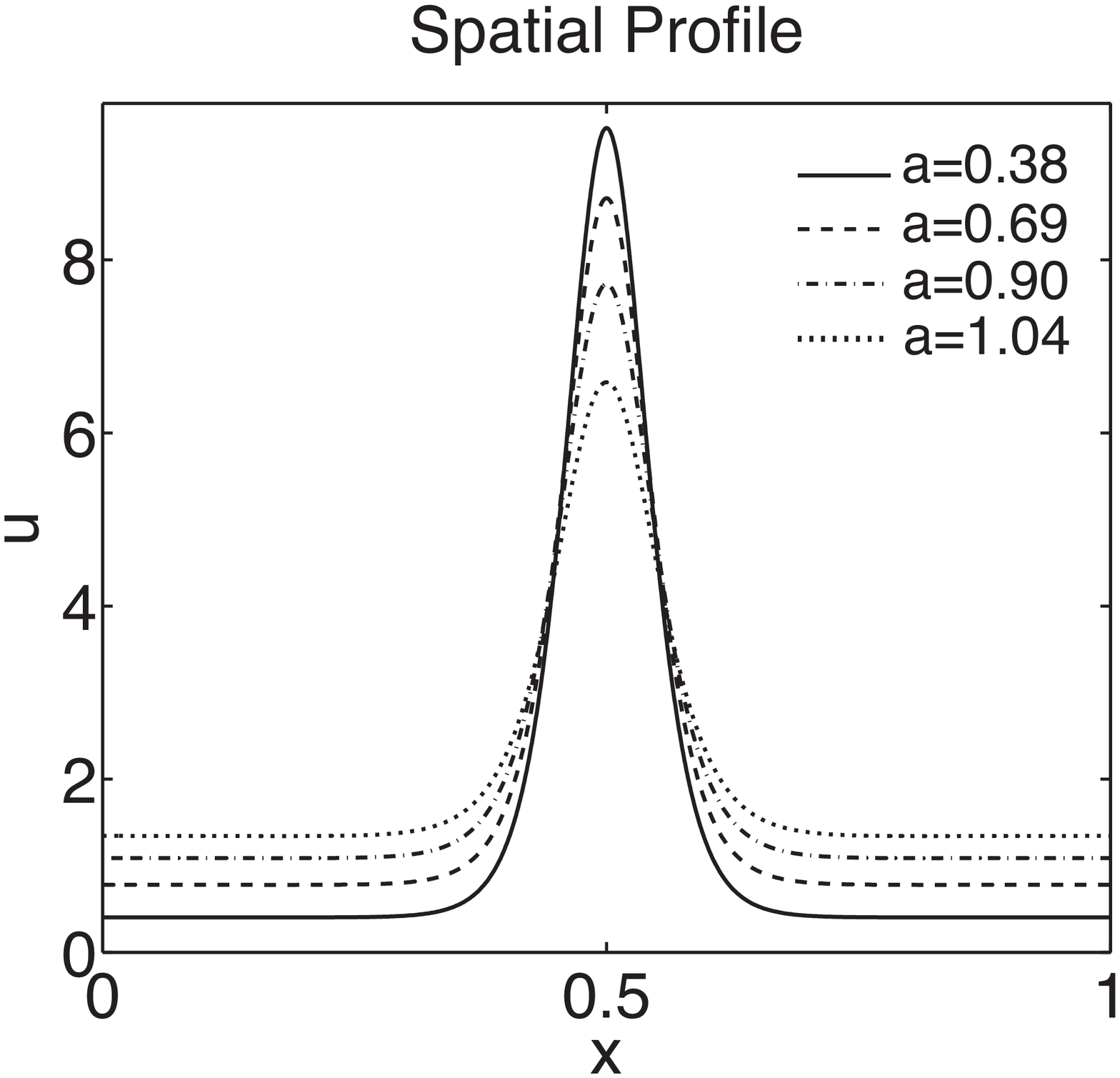}
}
\caption{Comparison of linear stability, local perturbation, and full PDE bifurcation analysis results for the Schnakenberg system \eqref{snack-sys}.  All diagrams are computed with $D=10$, $b=1$.  \textbf{Panel a) Local perturbation analysis results:} Global $u^s$ and local $u^{l1}$ solution branches of Eqn.~~\eqref{snack-LPA-sys} along with their stability are plotted as a function of $a$.  Two pattern forming regimes are predicted; I) linearly (Turing) unstable, and II) linearly stable where a sufficiently large perturbation induces patterning. \textbf{Panel b) Linear stability analysis results}  Maximum eigenvalue of $J_1$ \eqref{J_k} as a function of ``$a$'' for various values of $\epsilon$.  As $\epsilon \to 0$, the edge of the Turing region, marked with dots, approaches a limiting point near $a=1$, in agreement with LPA predictions in panel \emph{a}.  \textbf{Panel c) PDE bifurcation results:}  Bifurcation analysis of the full system of PDEs \eqref{snack-sys}.  The vertical axis describes the height (maximum - minimum) of a patterned solution.  As predicted by the LPA, stable patterned solutions exist both in the linearly unstable and stable regimes for sufficiently small $\epsilon$.  The location of the Turing bifurcation near $a=1$ agrees with panels \emph{a,b}.  Marked points represent points where the computed solution is plotted in panel \emph{d}.  \textbf{Panel d}) Example solutions of Eqn.~~\eqref{snack-sys} with $\epsilon=.025$.  }
\label{fig:Schnak}
\end{figure}

\noindent Based on the asymptotics above, it is expected that the initial behaviour of a local perturbation of the RDE's mimics the behaviour of the perturbation $u^l$ determined by this bifurcation analysis.   In region I ($a<1$ in Figure \ref{fig:Schnak}a), arbitrarily small perturbations of $u^l=u^s$ grow, predicting the HSS of Eqn.~~\eqref{snack-sys} is linearly unstable.  In region II, sufficiently large perturbations of $u^l$ are required to elicit a response for the LPA-ODE's, suggesting the HSS is linearly stable but large perturbations yield a response.

\vspace{.1in}
\noindent \textbf{Prediction 2:} In region II, as $a$ increases, increasingly large perturbations are required to initiate patterning.
\vspace{.1in}

\noindent In region II ($a>1$), the gap between the stable global and unstable local branches represents a response threshold for the LPA-ODE's: a perturbation of $u^l$ below the threshold decays back to the global HSS branch, a perturbation above it grows.  The dependence of this response threshold on the system parameter $a$ can be found by visual inspection of the LPA diagram.  For $a>1$, that threshold increases with $a$.  This threshold is precise only in the $\epsilon \to 0, D \to \infty$ limit, however the qualitative dependence on $a$ is expected to hold for the RDE system \eqref{snack-sys} with sufficiently extreme diffusivities.  

\vspace{.1in}
\noindent \textbf{Prediction 3:} The Turing bifurcation near $a=1$ is \textit{sub-critical} for sufficiently extreme diffusivities with large amplitude patterned states present on both sides of the bifurcation.
\vspace{.1in}

\noindent In dynamical systems theory, the terms \textit{sub-critical} and \textit{super-critical} are often used to describe the character of bifurcations such as Hopf or pitchfork.  Super-critical denotes a bifurcation that gives rise to a small amplitude response upon crossing it.  Sub-critical denotes one where the HSS loses stability immediately giving way to a large amplitude response.  In the latter case, responses can occur even outside of the unstable regime given a sufficient perturbation.  The Turing bifurcation near $a=1$ is predicted to be sub-critical with the unstable local branch $u^{l1}$ characterizing a threshold that shrinks to $0$ at the bifurcation, giving rise to instability of the HSS.

\vspace{.1in}
\noindent \textbf{Prediction 4:} Patterned solutions take the form of a spatially localized spike.
\vspace{.1in}

\noindent 
In region II, the LPA-ODE's exhibit blow up; a perturbation above the critical threshold grows to infinity.  When this perturbation becomes large, diffusion is expected to become important for the RDE's.  This will tend to oppose growth and smooth the solution.  It is reasonable to conjecture that at a particular height, reaction driven growth and diffusion driven suppression will balance leading to a large amplitude, spatially localized spike.  For future reference, it is expected that when $\epsilon$ decreases, decreasing the strength of diffusion, this spike would be come taller and more localized.

\subsubsection{Verification of predictions}

Figures \ref{fig:Schnak} (b,c,d) show results of linear stability, full PDE bifurcation, and numerical analyses for Eqn.~~\eqref{snack-sys}.  These results support the predictions above with a few caveats discussed at the end of this section.

\vspace{.1in}
\noindent \textbf{Verification of prediction 1}
\vspace{.1in}

\noindent Detection of a Turing instability for $a<1$ agrees well with LSA results presented in Figure \ref{fig:Schnak}b.  There, eigenvalues of the linearized Jacobian (i.e. Turing growth rates) are plotted as a function of the bifurcation parameter $a$ for four successively smaller values of $\epsilon$.  Dots on Figure \ref{fig:Schnak}b indicate LSA bifurcations; these bifurcation values are recorded for two different values of $D$ in Table \ref{table1}.  The location of these bifurcation values appears to converge to the predicted value of $a=1$.
This supports point 1 in  Section \ref{Important-points} that the LPA detects linear instabilities.

\begin{table} \label{table1}

\caption{Value of ``$a$'' at the edge of the Turing region for the Schnakenberg system \eqref{snack-sys} with $b=1$ for various values of $\epsilon$.  Column two: values drawn from the marked points in Figure \ref{fig:Schnak}b.  Column three: similar values for $D=1000$.  The final row is the value of the bifurcation predicted by the LPA.  This bifurcation approaches that predicted by the LPA as $\epsilon \to 0$, $D \to \infty$ in agreement with Corollary \ref{corollary-max-Tur}. }

\footnotesize \begin{center}
    \begin{tabular}{|c|c|c|}
        \hline
    $\epsilon$ & Turing ($D=10$) & Turing ($D=1000$) \\ \hline
    $0.1$   & $0.76$ & $0.82$ \\ 
    $0.05$  & $0.88$ & $0.95$ \\
    $0.025$ & $0.91$ & $0.98$ \\
    $0.01$  & $0.93$ & $0.99$ \\ \hline
    $LPA$  & $1$ & $1$  \\
        \hline
    \end{tabular}
\end{center}

\end{table}

The LPA predicts that in the linearly stable region II, sufficiently large perturbations elicit a patterning response.  To test this, both full PDE bifurcation analysis and asymptotics are employed.  Figure \ref{fig:Schnak}c shows results of numerical continuation (using Auto  \cite{Auto-citation}) of patterned solutions of the full PDE system  \eqref{snack-sys} with the vertical axis depicting the height (maximum - minimum) of the patterned solution.  The horizontal axis depicts the unpatterned HSS (maximum - minimum=0).  For sufficiently small values of $\epsilon$, the stable patterned solution extends into the $a>1$ region where the HSS is linearly stable, verifying the presence of stable patterned solutions in that regime.  Further, asymptotic results in Appendix \ref{Schnak-asym} show that in the $\epsilon \to 0, D \to \infty$ limit, patterned solutions exist for all values of $a$.  This lends support for point 2 in Section \ref{Important-points} that the LPA detects inherently non-linear patterning regimes.

\vspace{.1in}
\noindent \textbf{Verification of prediction 2}
\vspace{.1in}

\begin{figure}[htb] 
  \centerline{ 
\includegraphics[scale=.28]{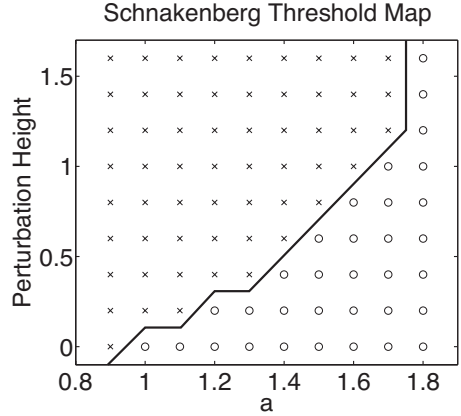}}
\caption{Verification of the qualitative relationship between ``$a$'' and the patterning threshold in the Schnakenberg system \eqref{snack-sys} using numerical simulation.  $\epsilon=.01$ and all other parameters are as in Figure \ref{fig:Schnak}.  Simulations were given a period of time to settle into a stable steady state (when present).  Perturbations of varying size were then applied to the middle $10 \%$ of the domain.  The vertical axis represents the size of the applied perturbation.  `x' indicates the perturbation grows resulting in a spike.  `o' indicates the perturbation decays back to the homogeneous state.  The patterning threshold increases with $a$ as indicated by the LPA.  To the left $a \approx 1$, the homogeneous state is unstable and to the right of $a \approx 1.7$, it is stable to all perturbations, in agreement with Figure \ref{fig:Schnak}c. }
\label{fig:SchnakThreshold}
\end{figure}

\noindent In Figure \ref{fig:SchnakThreshold} local perturbations of different height were applied to the HSS for multiple values of $a$ and the presence / absence of a pattern was recorded.  As $a$ increases, the perturbation size required to induce patterning increases as predicted by the LPA.  This supports point 4 in Section \ref{Important-points} that results of the LPA can be used to determine the qualitative dependence of response thresholds on parameters in non-linear patterning regimes.

\vspace{.1in}
\noindent \textbf{Verification of prediction 3}
\vspace{.1in}

\noindent  Figure \ref{fig:Schnak}c shows that as $\epsilon \to 0$, the nature of the Turing bifurcation changes from being super-critical to sub-critical.  For large $\epsilon$, small amplitude patterns emanate from the bifurcation.  For smaller $\epsilon$, the stable HSS gives way to large amplitude pattens immediately upon crossing the bifurcation.  Also for small $\epsilon$, an unstable patterned state is present outside the linearly unstable regime.  As the bifurcation is approached, this unstable state collapses onto the HSS changing its stability.  This is similar to the standard example of a sub-critical Hopf bifurcation where an unstable limit cycle colliding with a stable node results in a bifurcation.

It seems this association of a non-linear patterning regime with a sub-critical Turing bifurcation is somewhat common.  Both the substrate inhibition example and the chemotaxis example in Section \ref{GTPase-sec} presented later show sub-critical bifurcations.  Similarly, Rodrigues et al. \cite{Rodrigues-11} observed the presence of a stable heterogenous state adjacent to Turing parameter regimes for a discrete predator prey model.  They also showed the response threshold necessary to induce patterning in these regimes decreases as the Turing bifurcation is approached, consistent with these results.

\vspace{.1in}
\noindent \textbf{Verification of prediction 4}
\vspace{.1in}

\noindent Figure \ref{fig:Schnak}d and results in Appendix \ref{Schnak-asym} show that in both the linearly stable and unstable regimes a spike like solution forms.  Furthermore, these results show that as $\epsilon$ is decreased (reducing the opposing effect of diffusion),  the spike height increases as expected.  Thus the LPA inferences about the long term evolution of patterns are confirmed in this example, supporting point 3 in Section \ref{Important-points}.

\subsubsection{Caveats of the local perturbation results}

There are discrepancies between the LPA predictions and results of LSA and full PDE bifurcation analyses.  First, the location of the predicted bifurcation at $a=1$ is not precise.  Both LSA and PDE bifurcation results show the value of the actual Turing bifurcation depends on $\epsilon$ and $D$.  Though this does appear to converge to the predicted $a=1$ in the proper limit.  This type of approximation error will be present in any LPA application and will be discussed in more detail in Section \ref{sec:LPA-LSA}.

Second, the LPA predicts patterned solutions will form for all $a>1$.  Full PDE bifurcation results (Figure \ref{fig:Schnak}c) in contrast show the patterned state  is annihilated in a saddle node (or fold) bifurcation at a finite value of $a=a^*(b,\epsilon,D)$.  The location of this bifurcation does increase as $\epsilon \to 0$ and asymptotic results in Appendix \ref{Schnak-asym} show the presence of a patterned solution for all values of $a$.  So in the $\epsilon \to 0, D \to \infty$ limit, $a^* \to \infty$, in agreement with LPA results.

\subsection{The local perturbation analysis of a substrate inhibition model}\label{SubInh-results}
I now apply the LPA to a substrate inhibition model to demonstrate a different set of results and interpretations obtained with the same method.
This model \cite{Thomas-79}
\begin{subequations}\label{SI-sys}
\begin{align}\label{SI-sys-u}
u_t(x,t)&=a-u-\frac{\rho u v}{1+u+Ku^2}+\epsilon^2 \triangle u = f(u,v) +\epsilon^2 \triangle u, \\
v_t(x,t)&=\alpha (b-v)-\frac{\rho u v}{1+u+Ku^2}+ D \triangle v =g(u,v)+D \triangle v ,\label{SI-sys-v}
\end{align}
\end{subequations}
describes two co-substrates that are constantly generated, decay linearly, and are used up in an enzymatic reaction.  The non-linear term is indicative of multiple substrate molecules $u$ binding to a single enzyme rendering it inert for further interaction with the remaining co-substrate $v$, thus the term substrate inhibition.  See \cite{Murray-82} for LSA results for this system.

\begin{figure}[htb!]
\psfrag{Du}{$D_u$}
\subfigure[]{
\centering
\includegraphics[scale=.23]{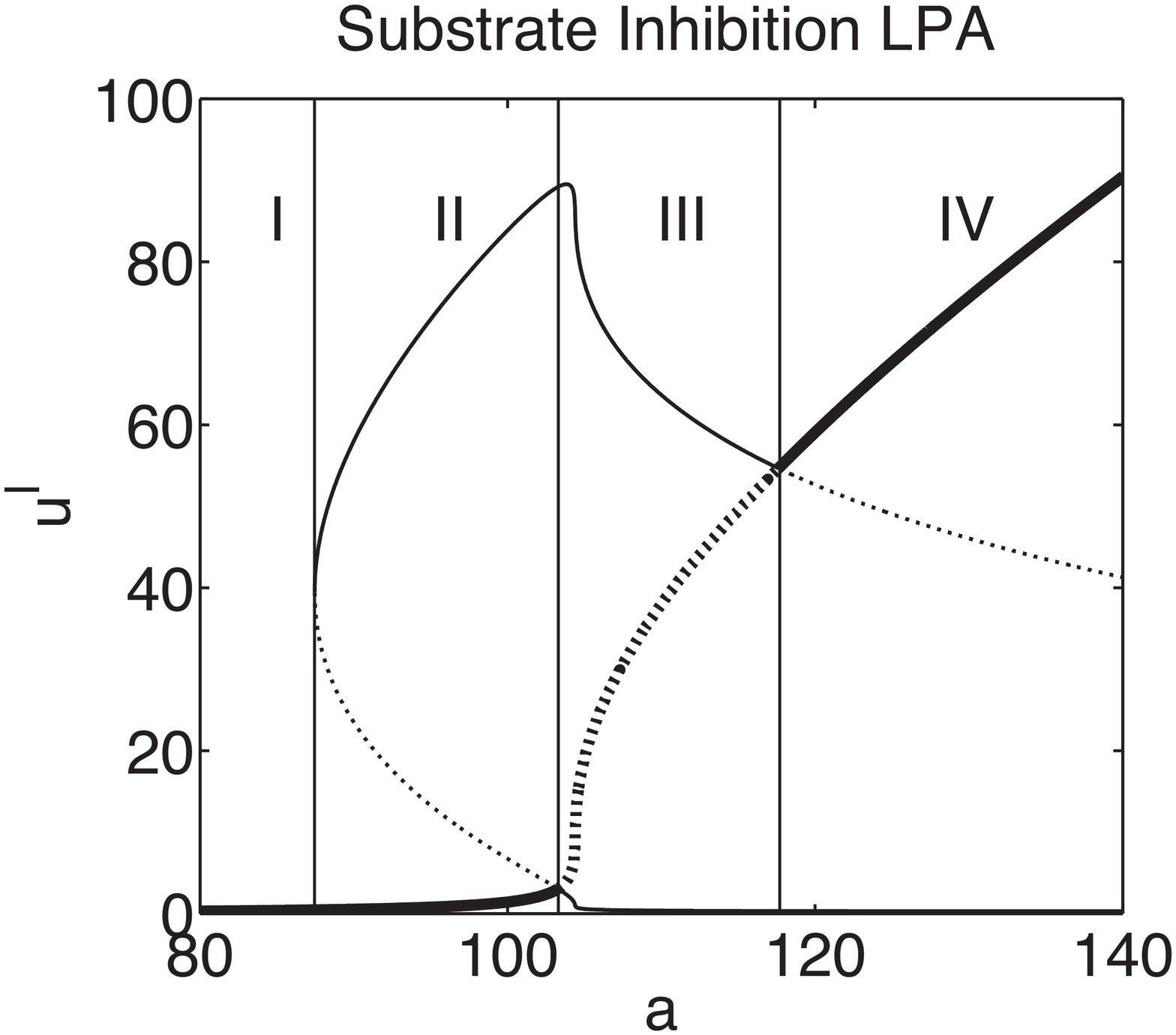}
\label{SI-a}
}
\subfigure[]{
\centering
\includegraphics[scale=.23]{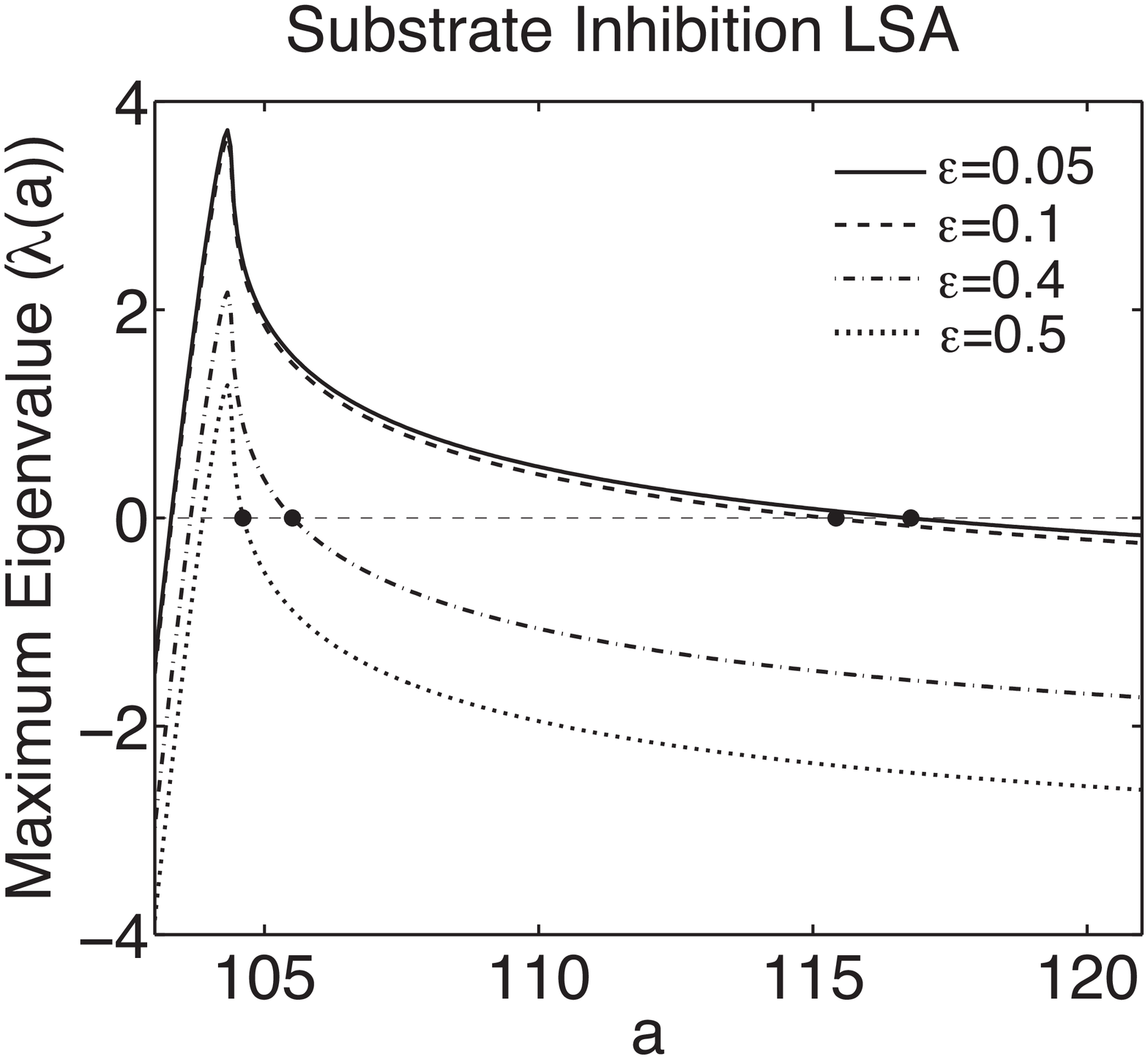}
\label{SI-b}
}
\subfigure[]{
\centering
\includegraphics[scale=.23]{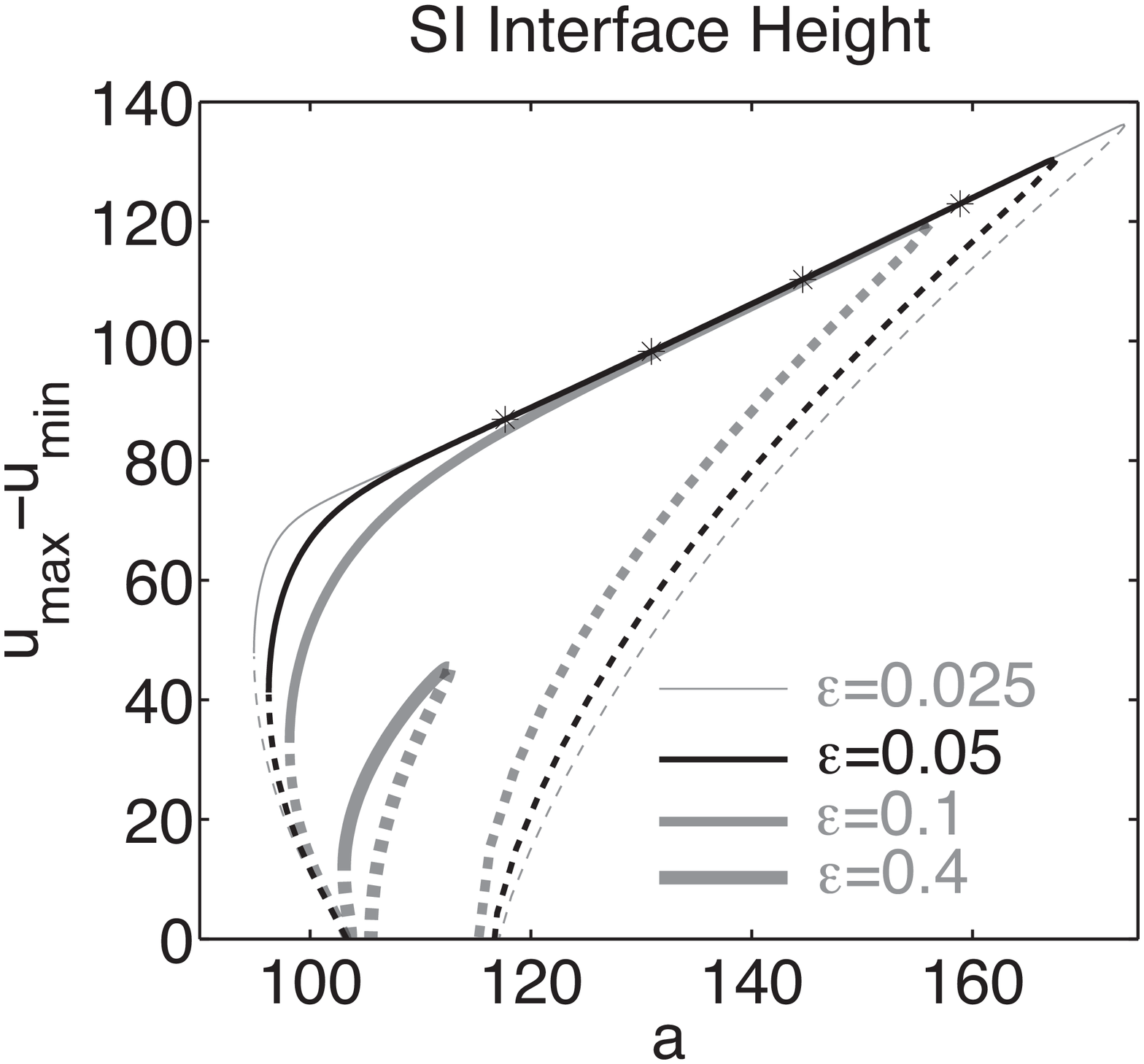}
\label{SI-c}
}
\subfigure[]{
\centering
\includegraphics[scale=.23]{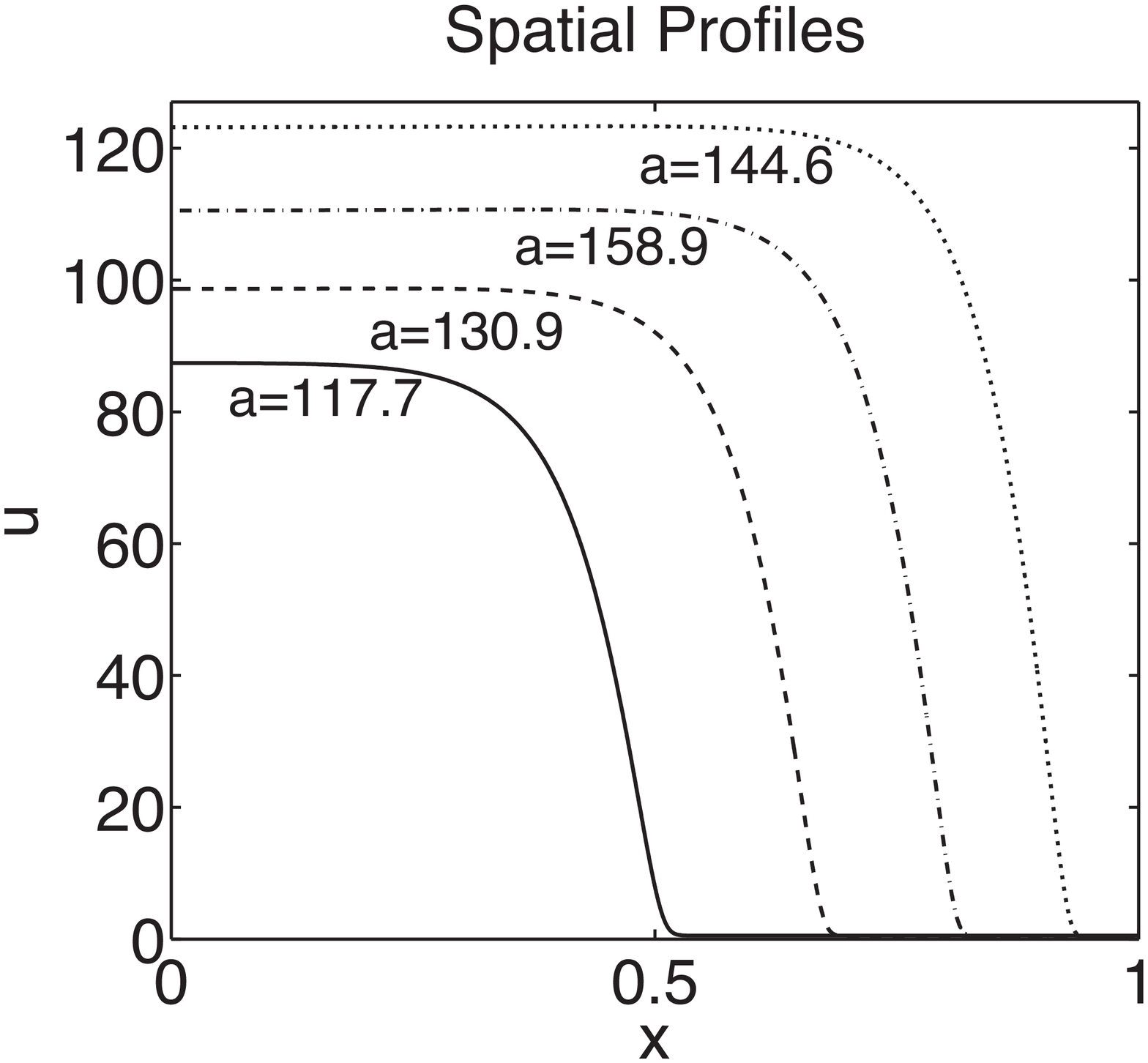}
\label{SI-d}
}
\caption{\textbf{Comparison of linear stability, local perturbation, and full PDE bifurcation analysis results for the substrate inhibition model Eqn.~~\eqref{SI-sys}}.  All diagrams are computed with $D=10$, $\rho=13$, $K=0.125$, $\alpha=1.5$, and $b=80$ and all conventions are as in Figure \ref{fig:Schnak}.  \textbf{Panel a) Local perturbation analysis results:}  Three regimes of behaviour are predicted: No patterning (I), linearly (Turing) unstable (III), and linearly stable where a sufficiently large perturbation yields a response (II, IV).  \textbf{Panel b) Linear stability analysis results:} For each $\epsilon$ two Turing bifurcations are seen.  Dots mark the location of the right Turing bifurcation for different values of $\epsilon$.  As $\epsilon \to 0$, the edges of the Turing region approach that predicted by the local perturbation analysis results in panel \emph{a}.  \textbf{Panel c) PDE Bifurcation Results:}  In agreement with the LPA results, patterned solutions are found both inside and outside the linearly unstable regime.  \textbf{Panel d)} Example solutions drawn from starred points for $\epsilon=.05$ in panel \emph{c}.  As predicted by the LPA results, these solutions show a stable interface separating high / low regions of $u$. }

\label{fig:SubInh}
\end{figure}

In the previous example, it was possible to analytically compute the various solutions of the LPA system of ODE's.  This will not generally be the case, but it is possible to find and track the various LPA solution branches efficiently using standard ODE bifurcation techniques.  Figure \ref{fig:SubInh} mirrors Figure \ref{fig:Schnak} with results in \emph{panel b,c,d} verifying predictions inferred from LPA results in \emph{panel a}.
In this example, solution branches of the LPA-ODE's along with there stability are computed with the numerical continuation software package Matcont \cite{Matcont-03}.  LPA predictions regarding the presence of linearly unstable (region III) and non-linear (regions II, IV) pattering regimes are again verified by LSA and full PDE bifurcation analyses.  

There are a few important contrasts between these LPA results and those for the Schnakenberg model.  First, in region IV (Figure \ref{fig:SubInh}a), a negative valued perturbation of the HSS is predicted to induce patterning.  This along with the general dependence of response thresholds in regions II, IV were verified numerically (results not presented).  Second, the LPA results suggest the resulting solution will take the form of a stable interface separating high and low regions for $u$.

Consider region II where this system has two local branches, one unstable and the other stable.  In this case, a perturbation of $u^l$ above the unstable local branch will be attracted to the stable local branch.  On the slow diffusion timescale, boundary layers between the high and low $u$ states on $R^{l,g}$ respectively will move.  The lack of a higher HSS suggests the high state cannot encompass the entire domain and one of two things will happen: 1) the boundary layers will stop / stall somewhere in the interior of the domain leaving stable interfaces separating regions of high and low $u$ levels, or 2) the solution will eventually collapse back to the unique HSS.  Figure \ref{SI-d} showing stable long term solutions for multiple values of $a$ verifies the former is the case here.  These results again support the suppositions in Section \ref{Important-points}.

\section{Detection of linear instabilities by the LPA}\label{sec:LPA-LSA}

The previous section demonstrated points 1-4 in Section \ref{Important-points}.  Points 2-4 relate to non-linear patterning regimes and will be difficult to show in generality.  The supposition in point 1 does however hold in generality.  In this section I show the LPA detects linear instabilities with increasing accuracy as diffusivities become extreme.  Asymptotic analysis shows eigenvalues of the linearized RDE's converge to eigenvalues of the linearized LPA-ODE's.  Thus positive eigenvalues of the LPA system (indicating linear instability) correspond to positive Turing growth rates and vica versa.



Since eigenvalues of the linearized RDE's \eqref{LPA-1} determine linear stability properties, they will be characterized in the $\epsilon \to 0, D \to \infty$ limit.  Linearizing about the HSS $(u^s,v^s)$ with respect to periodic perturbations of the form $\exp (i k x)$ leads to the Jacobain
\begin{equation}\label{J_k}
J_k = \begin{bmatrix} f_u (u^s,v^s;p) - k^2 \epsilon^2 I & f_v(u^s,v^s;p) \\ g_u(u^s,v^s;p) & g_v(u^s,v^s;p) - k^2 D I \end{bmatrix} .
\end{equation}
Recall that $f: \mathbb{R}^M \times \mathbb{R}^N \to \mathbb{R}^M, g: \mathbb{R}^M \times \mathbb{R}^N \to \mathbb{R}^N$ and denote eigenvalues of this matrix as $\{ \lambda^k_i(\epsilon,D,p) \}_{i=1:(M+N)}$.  Assume these are in decreasing order according to their real part so that $\lambda^k_1$ has the largest real part and determines stability.

Stability of the HSS branch of the LPA ODE's $(u^g,v^g,u^l)=(u^s,v^s,u^s)$  is determined by eigenvalues of the LPA Jacobain
\begin{equation}\label{LP-J}
J_{LP}=
\begin{bmatrix}
f_{u}(u^s,v^s;p) & f_{v}(u^s,v^s;p) & 0 \\
g_{u}(u^s,v^s;p) & g_{v}(u^s,v^s;p) & 0 \\
0 & f_{v}(u^s,v^s;p) & f_{u}(u^s,v^s;p) 
\end{bmatrix} 
.
\end{equation}
This is a block lower triangular matrix with the upper left block being precisely $J_0$ \eqref{J_k}.  Thus the eigenvalues $\{ \lambda^0_i \}_i$ of $J_0$ are also eigenvalues of $J_{LP}$.  These describe stability of the HSS of Eqn.~~\eqref{LPA-1} with respect to spatially homogeneous perturbations and determine stability of the well mixed system (i.e. with $D_u = 0 = D_v$).   So the LPA detects any instabilities of the well mixed system.  These will be referred  to as the well mixed eigenvalues of $J_{LP}$.   

The remaining eigenvalues of $J_{LP}$ are simply eigenvalues of $f_u(u^s,v^s)$.  Denote these as $\{ \lambda^{LP}_j(p) \}_{j=1:M}$ and assume they are ordered according to decreasing real part so that $\lambda^{LP}_1$ has the largest real part and determines stability.  Then the following asymptotic result relating $\{ \lambda^{LP} \}$ and $\{ \lambda^k \}$ holds.
\vspace{.1in}
\begin{thm}\label{Thrm-EigVal}
Assume $\epsilon^2 \ll 1 \ll D$, $\nabla f$, $\nabla g$ are $O(1)$ with respect to $\epsilon$ and $D$, and fix a wave number $k>0$.  Further assume that $f_u (u^s,v^s), g_v(u^s,v^s)$ are diagonalizable.  Then:
\begin{enumerate}

\item For each $i=1:M$,  $\lambda^k_i=\lambda^{LP}_i -k^2 \epsilon^2 +c(D)$ where $c(D) \to 0$ as $D \to \infty$.

\item For each $i=M+1:M+N$, $\Re (\lambda^k_i) = O^{-}(D)$ where $O^{-}()$ signifies a negative valued quantity of that order.

\end{enumerate}
\end{thm}
\vspace{.1in}

\noindent For proof of this result, see Appendix \ref{Thrm-EigVal-proof}.  A direct consequence of this is that as $\epsilon \to 0, D \to \infty$, $\lambda_1^k(\epsilon,D,p) \to \lambda_1^{LP}(p)$ independent of $k$.  So in this limit, linear instability of the HSS branch of the LPA-ODE's corresponds directly to diffusion driven (i.e. Turing) instability for the RDE's.  Thus point 1 in Section \ref{Important-points} holds for general systems of the form \eqref{LPA-1} when $f_u (u^s,v^s),g_v (u^s,v^s)$ are diagonalizable.

\subsection{LPA bifurcations locate the edge of ``limiting'' linearly unstable parameter regimes}

Recall from the LSA results in the previous examples that the location of a Turing bifurcation of the RDE's appears to converge to a limiting point as $\epsilon \to 0, D \to \infty$.  This is in line with Murray's \cite{Murray-82} observation that a linearly unstable regime of parameter space converges to a ``limiting'' unstable regime in this limit.  As indicated by the comparison of the locations of LSA bifurcations and those predicted by the LPA, the LPA precisely locates the edge of these ``limiting'' unstable regimes.  This follows from the following corollary of Theorem \ref{Thrm-EigVal}.

\vspace{.1in}
\begin{cor}\label{corollary-max-Tur}
Consider a particular set of parameters $p$.  If the global branch of the LPA-ODE's is linearly unstable ( i.e. $\Re (\lambda_1^{LP}(p) )>0$), then the HSS of the RDE's is linearly unstable for sufficiently extreme values of $\epsilon$ and $D$ (i.e. $\Re (\lambda_1^{k} (\epsilon,D,p))>0$ for some $k>0$).  Furthermore, if the global branch is linearly stable in the LPA sense, the HSS is linearly stable for sufficiently extreme values of $\epsilon, D$ as well.

\end{cor}
\vspace{.1in}




\section{Applications of the LPA to more complex systems} \label{GTPase-sec}

One of the primary benefits of the Local Perturbation Analysis is the relative ease with which it can be applied to more complex systems, common in biological applications.  Existing methods become either difficult to implement or intractable in such cases.  However, with the help of ODE analysis software packages such as Auto \cite{Auto-citation} and Matcont \cite{Matcont-03}, this method scales well to larger systems.  The following example demonstrates an application of the LPA to a system involving 9 RDE's.


\subsection{Chemotactic polarization example}

Much effort has been devoted to understanding the process by which cells, ranging from white blood cells to cancer cells, move up chemical gradients.  Reorganization of regulatory molecules, primarily GTPases and phosphoinositides, is known to be a precursor to such motion.  In response to an applied chemical gradient, these molecules self organize to form a polar state where some localize in the cell ``front'' (Cdc42, Rac, PI3K, and PIP3) and others in the ``rear'' (Rho, PTEN).  Front related molecules generate protrusion, rear related molecules generate contraction, and their combined activity leads to directed motion.

Each of the three GTPases (Cdc42 $C$, Rac $R$, and Rho $\rho$) effectively has two forms, membrane bound and cytosolic with only the membrane bound form in an active state.  Over the timescale of polarization events, the amount of each GTPase is conserved with diffusion and cross talk mediated cycling between the two states leading to segregation of active forms.  These cross talk interactions and the influence of phosphoinositide feedback are the focus of this discussion.

While these regulators are conserved across a wide range of eukaryotic cells, the cross talk interactions between them is not.  This variation has led to extensive experimental work aimed at dissecting these interactions in different cell types and numerous models (reviewed in \cite{Jilkine-11}) aimed at understanding their results.  Here I describe and analyze a variant of a model \cite{Holmes12a,Lin-12} motivated by work on HeLa cell polarization.

I investigate a structural perturbation of that model, introducing mutual antagonism between Rac and Rho, known to be present in numerous cell types \cite{Sanders-99,Caron-03,vanLeeuwen-97}.  A schematic diagram of this model is in Figure \ref{fig:GTPase2}a.  The dashed interaction, Rho mediated inhibition of Rac, is the structural addition differentiating this model from that in \cite{Holmes12a,Lin-12}.  Model equations encoding these interactions are found in Eqs.~\eqref{equ:GTPase} with a description of parameters and their values in Table \ref{table:params}.  See Appendix \ref{GTPase-equ} for a brief description of this model and \cite{Holmes12a,Lin-12} for more extensive discussion of the original network and its parameters.

What effect does the addition of Rho mediated inhibition of Rac have on the behavior of this network?  To investigate this, a non-dimensional parameter ($f_2$) modulating the strength of this inhibitory interaction is introduced.  When $f_2=0$, no inhibition is present and the original network is recovered.  When $f_2$ increases, the strength of the inhibition increases.  LPA and numerical simulation results in Figure \ref{fig:GTPase2} show the effect of increasing the strength of this feedback.

\subsection{LPA and numerical simulation results}

\begin{figure}[htb!]
\psfrag{Du}{$D_u$}

\subfigure[]{
\centering
\includegraphics[scale=.39]{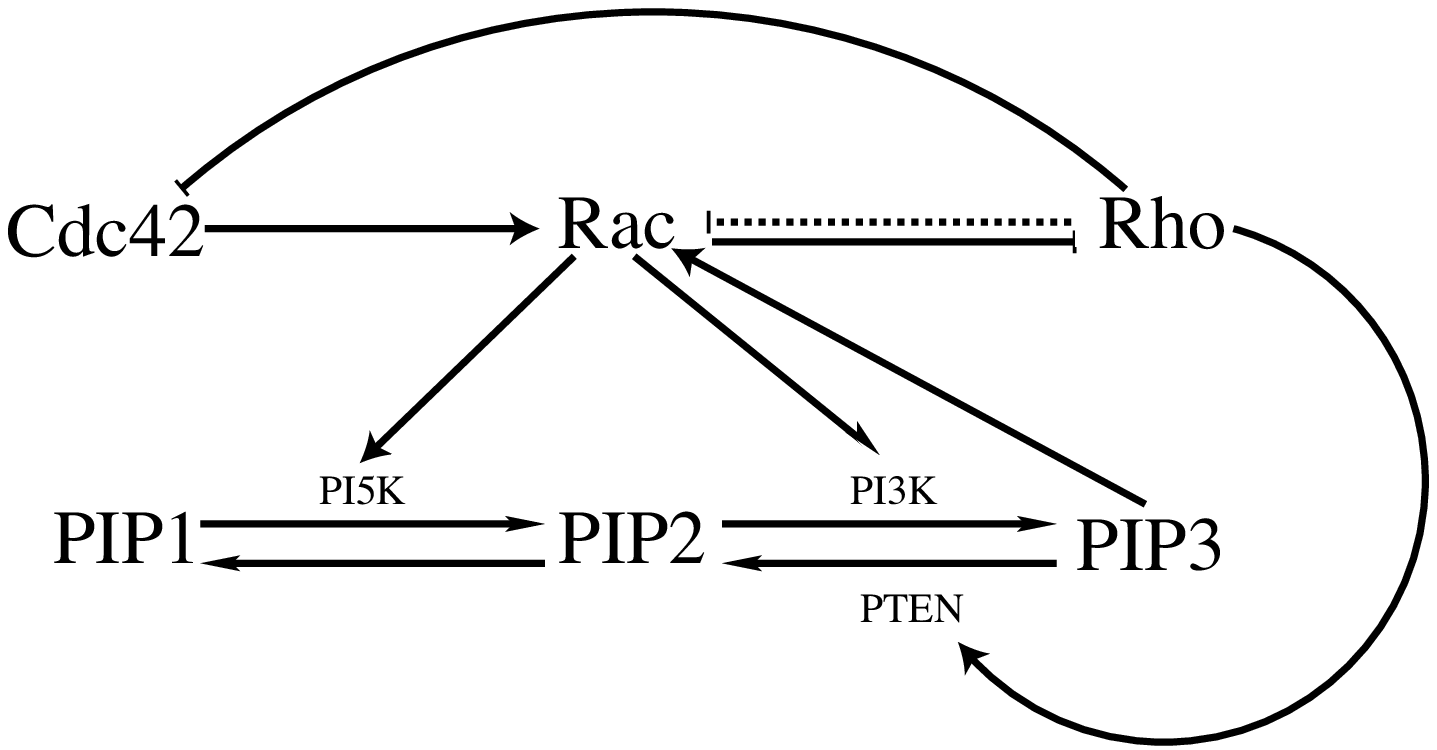}
}
\subfigure[]{
\centering
\includegraphics[scale=.22]{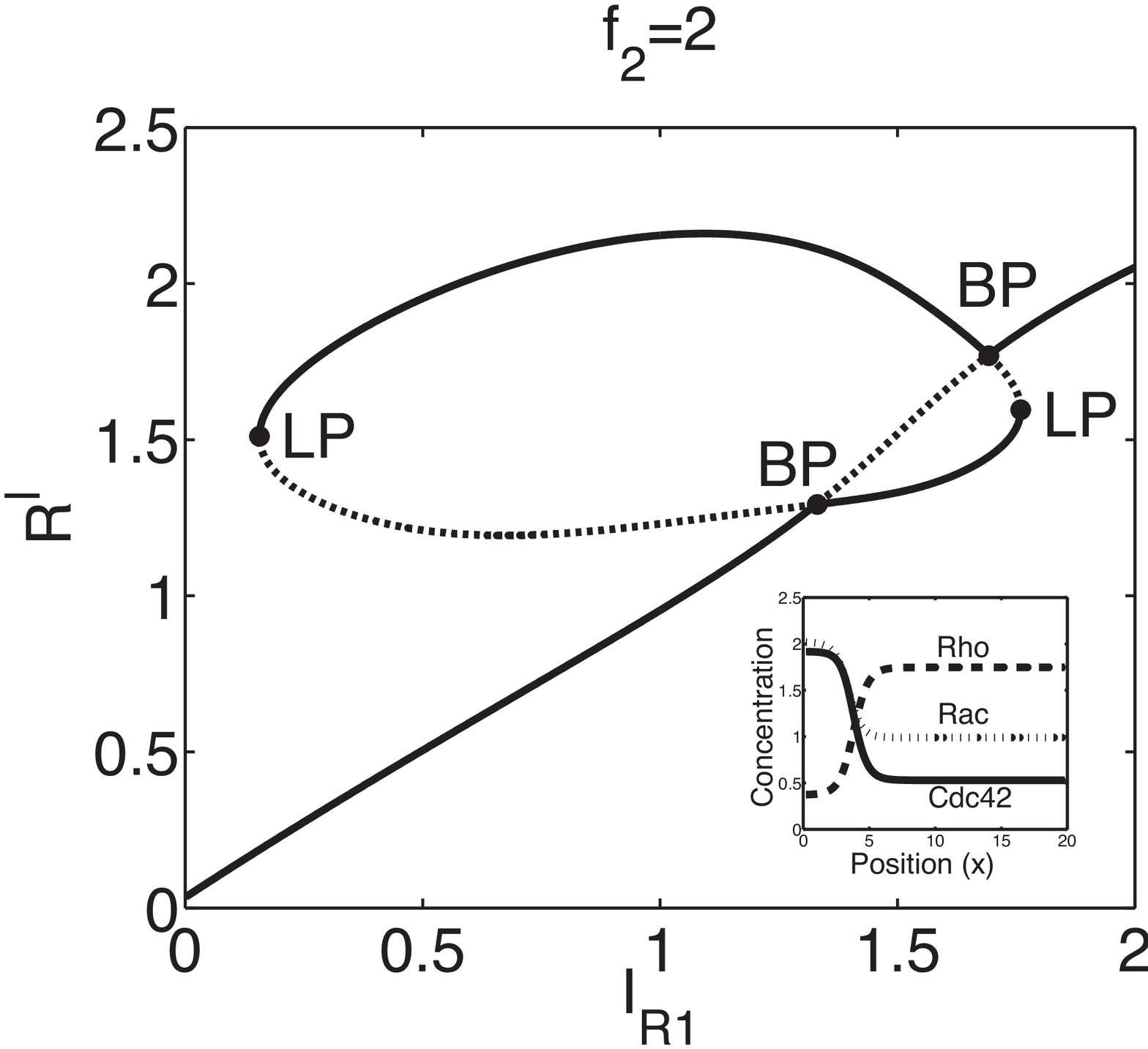}
}
\subfigure[]{
\centering
\includegraphics[scale=.217]{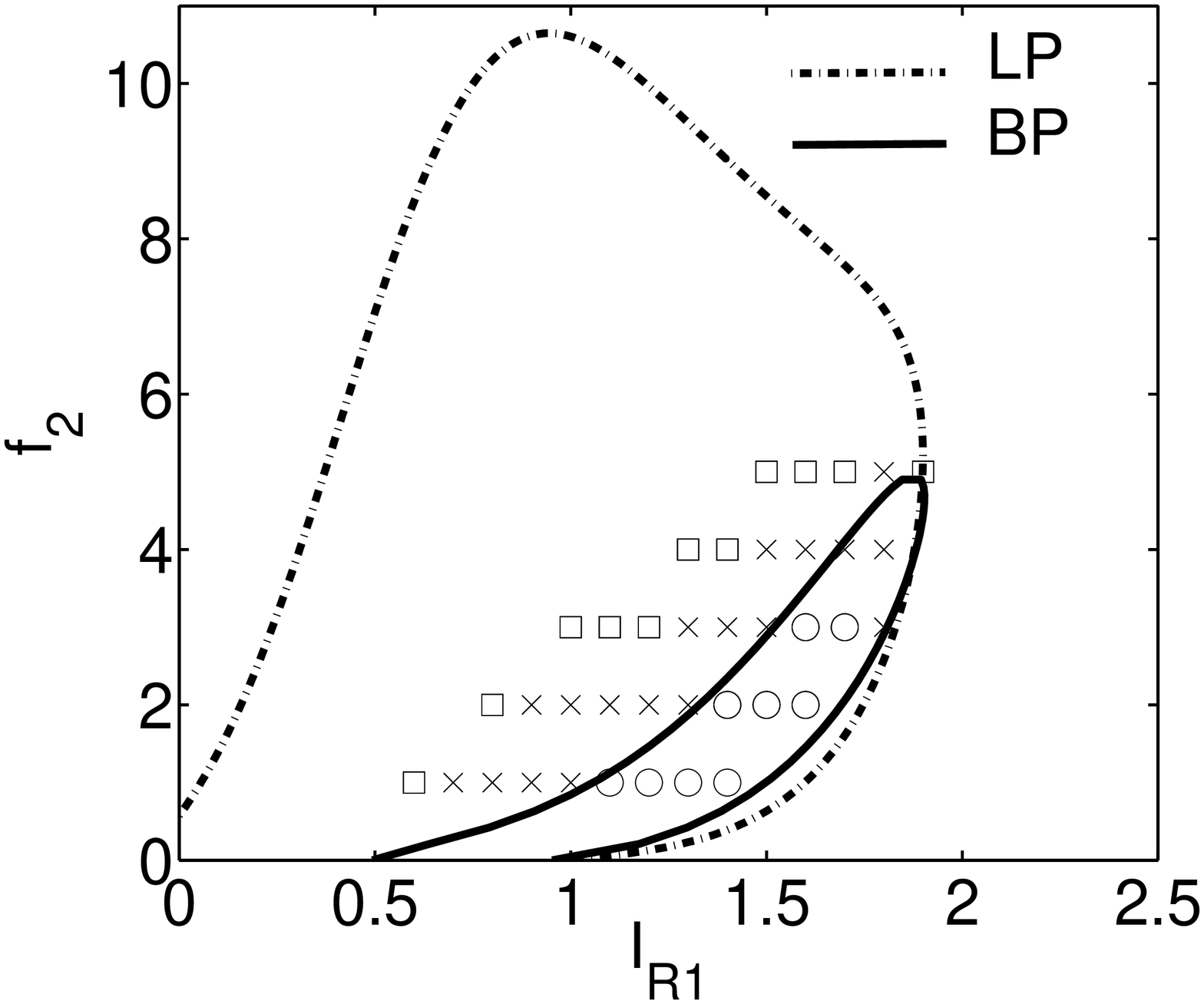}
}
\subfigure[]{
\centering
\includegraphics[scale=.23]{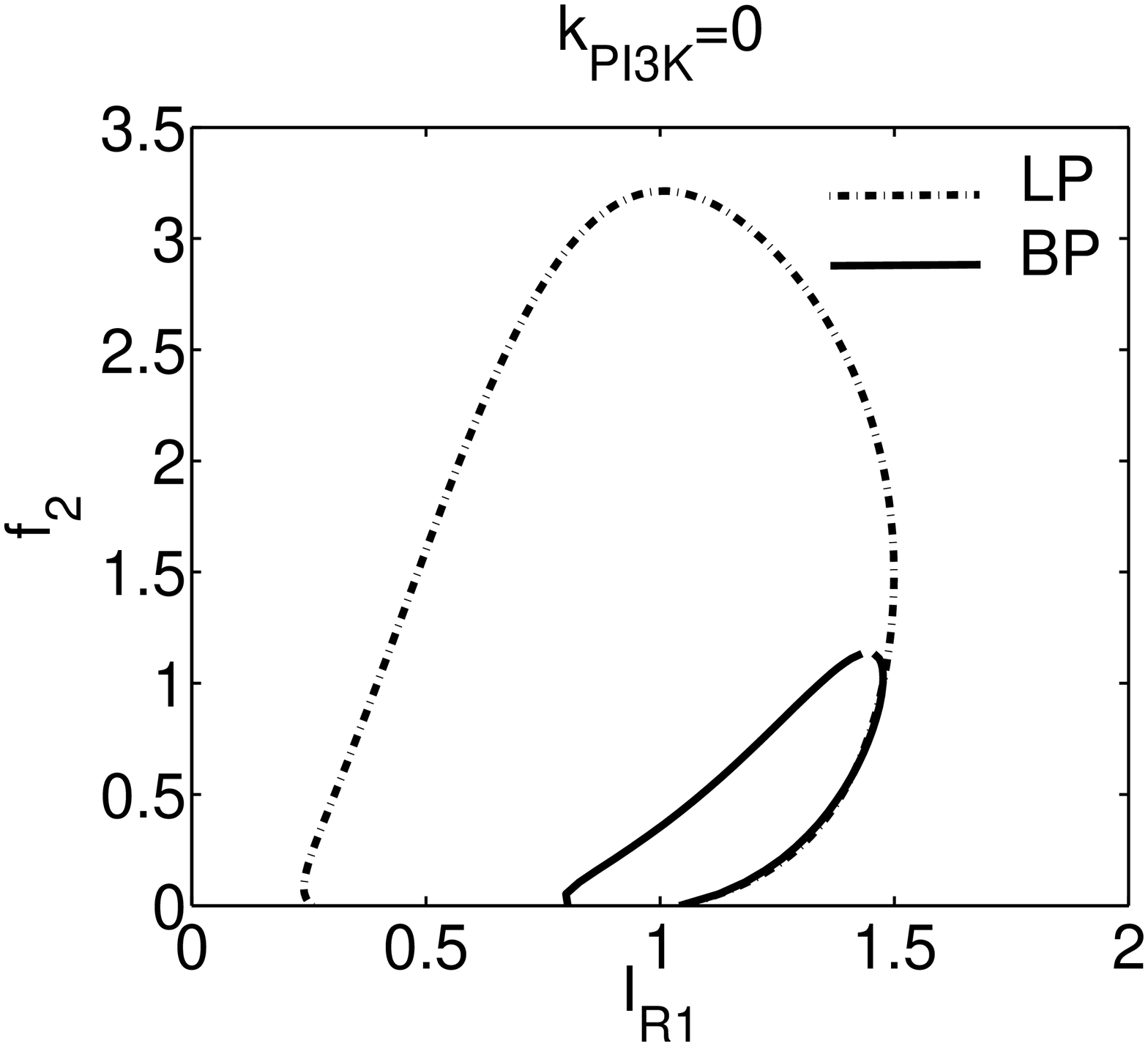}
}
\label{fig:GTPase2}
\caption{\textbf{Panel a} Schematic diagram of interactions between three GTPases (Cdc42, Rac, and Rho) along with phosphoinositides (PIP1, PIP2, PIP3) represented by Eqs.~\eqref{equ:GTPase}.  An arrow represents activation of one protein by another, a bar represents inactivation.  \textbf{Panel b:} Local perturbation analysis of this system with $f_2=2$ and $I_{R1}$ the bifurcation parameter.    The vertical axis is the value of the local form $R^l$ of Rac.  The monotonic branch is the global branch, the loop is the local branch   Four regions are seen, ordered from left to right: stable, linearly stable where a sufficiently large perturbation yields a response, linearly unstable, and stable.  There is a small nonlinear patterning regime regime between $I_{R1} \sim 1.6, 1.7$, but this regime appears to only be present for extreme diffusions beyond those considered here.  \textbf{Inset:} Numerical simulation results at $f_2=2$, $I_{R1}=1.1$.  \textbf{Panel c:}  Two parameter continuation of the branch points (BP), indicating the edge of a linearly unstable regime, and fold bifurcations (LP) of the local branch.  Markers represent simulation results of the full RDE system.  Circles indicate a Turing instability, `x' indicates perturbation induced patterning, and a square indicates a parameter set which is stable for $D_m=.1$ but where perturbation driven patterning results with $D_m=.01$.  \textbf{Panel d:} The same as Panel c with PI3K knockdown removing the feedback of PIP3 $\to$ Rac.  Similar parameter space structure is present though substantially compressed.  Values for parameters not explicitly mentioned are in Table \ref{table:params}. }
\end{figure}

Figure \ref{fig:GTPase2}a shows the results of a LPA of this model with moderate feedback, $f_2=2$.  The LPA was performed assuming  membrane bound GTPases are slow diffusing ($D_{m}=0.1 \mu m^2/s$), and cytosolic GTPases ($D_{c}=50 \mu m^2/s$) are fast.  For reference, cell sizes considered are on the order of $10-20 \mu m$.  Phosphoinositide diffusion lies between the fast and slow regimes.  However, as in \cite{Holmes12a}, LPA results are similar with it chosen as either fast or slow.  In Figure \ref{fig:GTPase2} they are taken to be slow variables.  The LPA reduction in this case leads to a system of 15 ODE's for 6 local variables ($C^l,R^l,\rho^l,P_1^l,P_2^l,P_3^l$) and 9 global variables ($C^g,R^g,\rho^g,C_c^g,R_c^g,\rho_c^g,P_1^g,P_2^g,P_3^g$).

The bifurcation parameter $I_{R1}$ represents a basal Rac activation rate.  Its variation could result from either population heterogeneity or external stimulation of Rac as in \cite{Lin-12}.  In Figure \ref{fig:GTPase2}b, three regimes of behavior are found at different activation levels.  For both low and high levels of basal activation, no response due to either instability or an applied stimulus can occur.  For increasing levels of activation, a regime where sufficiently large perturbations yield a response is found.  Again, this regime terminates in a sub-critical Turing bifurcation as the response threshold shrinks to zero.  This suggests increasing the spatially uniform activation rate increases the sensitivity of a cell to heterogeneous stimuli, experimentally supported in  \cite{Lin-12}.

At yet higher values of $I_{R1}$, a second narrow linearly stable patterning regime is found.  Numerical simulations verify the presence of all but this narrow regime, which likely requires more extreme diffusivities to be observed numerically.  LPA results again suggest solutions will evolve to a polarized profile with a transition layer separating regions of homogeneous activity levels.  This was also verified numerically with an indicative steady state solution shown in the inset ($I_{R1}=1.1$).

Now consider the effect of increasing / decreasing the strength of Rho $\dashv$ Rac feedback.  Labeled branch points (BP) indicate the approximate boarder of the unstable regime.  Fold bifurcations (LP) mark the boarder of the non-linear patterning regimes.  Standard two parameter continuation techniques are applied to follow these bifurcations as $f_2$ is varied, Figure \ref{fig:GTPase2}c.  At low values, both regimes persist.  As $f_2$ increases, the branch points collapse at $f_2 \sim 5$ and the linearly unstable regime between them is lost.  For higher values of $f_2$, the two fold bifurcations of the local branch persist suggesting the continued presence of a linearly stable patterning regime between them.

Marked points on Figure \ref{fig:GTPase2}c indicate parameter values where numerical simulation of the full RDE system was performed.  Circles indicate a parameter set where small noise (machine noise) induces a response.  Points marked $\times$ indicate sufficiently large perturbations are required for a response.  Beyond these points, no patterning was detected numerically at the base diffusion values.  When $D_m=.01 \mu m^2 / sec$, parameter regimes expand with squares marking additional parameter sets where sufficiently large perturbations yield a response.

The linearly unstable regime for the RDE's is confined to that predicted by the LPA and as expected, stimulus induced patterning is present to the left but not to the right of that regime.  While the location of patterning regimes in parameter space agree well with predictions, the expanse of these regimes is substantially smaller than predicted, particularly for the non-linear patterning regime.  However, as diffusivities are driven to yet further extremes, these regimes do expand further (results not presented).

I consider one final perturbation of this network, PI3K knockout which is accomplished here by setting $k_{PI3K=0}$.  This has the effect of removing feedback between GTPase's and phosphoinositide's.  The same analysis above was performed with results in Figure \ref{fig:GTPase2}d.  Similar parameter space structure exists with linearly unstable and stable patterning regimes.  In this case however they are substantially compressed in parameter space.  So while PI3K / PIP3 mediated feedback is not necessary for polarization, it does make it more robust in a parametric sense.  This is consistent with observations \cite{Ferguson-06,Lin-12} showing PI3K / PIP3 localization is not necessary for efficient chemotaxis, but its knockout substantially reduces the fraction of cells that do chemotax.

\section{Limitations of the LPA approximation}\label{sec:limitations}

Here I stress the limitations of the Local Perturbation Analysis.  The purpose of the LPA is not to approximate a solution of Eqn.~~\eqref{LPA-1}, only the initial response of a HSS to a localized perturbation \eqref{IC}, i.e. growth or decay.  Slow diffusion timescale and boundary layer effects are not considered  and the LPA-ODE's \eqref{LP-sys} only describe the evolution of the perturbation on short to intermediate timescales.  These effects can become important and the leading order approximation above can fail for one of two basic reasons.  
\begin{itemize}
\item The perturbation becomes large.  This would cause a number of effects.  First, if $g$ is unbounded, the correction term in Eqn.~~\eqref{vg-correction} could become $O(1)$ and affect leading order dynamics.  Second, if $f,g$ are unbounded, neglected Taylor expansion terms of the form $\epsilon f_u$, $\epsilon f_v$, $\epsilon g_u$, or $\epsilon g_v$ can become $O(1)$, influencing the leading order dynamics.  Third, boundary layer effects could influence the dynamics of $(U,V)$ on $R^l$.

\item The perturbation spreads in space.  This would again cause the area of the perturbed region to become $O(1)$, causing the correction term in Eqn.~~\eqref{vg-correction} to affect leading order dynamics.  These effects would however occur on the slow diffusion timescale which is not considered here.
\end{itemize}
In either case, these effects only become important after the perturbation has grown in amplitude, constituting a response.

As a result of neglecting these higher order effects, LPA predictions are asymptotic in nature and require sufficiently distinct diffusivities to be valid.  Therefore the predicted location of bifurcations between parameter regimes only approximate the location of actual bifurcations, with this approximation improving as diffusivities become more extreme.  In some cases, bifurcations of the RDE system with finite diffusivities (such as the saddle node bifurcation in the Schnakenberg case) are not captured at all by the LPA.  

Further, all diffusion related information is lost.  Therefore length scale information cannot be obtained, non-linear phenomena such as peak splitting \cite{Kolokolnikov-05,Wei-05,Barrass-06} will not be found, and pattern selection or abbarent effects from domain growth for example \cite{Crampin-99, Barrass-06} cannot be discussed.  For these reasons, the LPA should be viewed primarily as an efficient, scalable non-linear stability technique that can be used to inform further analysis.

\section{Discussion}
  
A new non-linear bifurcation technique for systems of reaction diffusion equations with large diffusion disparities \eqref{LPA-1} was developed and demonstrated.  This ``Local Perturbation Analysis'' (LPA) determines the response of a HSS of a system of reaction diffusion equations to a spatially localized,  large amplitude perturbation.  Under proper asymptotic assumptions about the diffusivities $\epsilon,D$ and the form of the perturbation, its evolution can be approximated to leading order by a collection of ODE's describing the perturbation (local variables) and the broader domain (global variables).

A bifurcation analysis of this collection of LPA-ODE's reveals two types of solution branches : 1) a ``global'' branch of solutions representing HSS solutions of the RDE's, and 2) ``local'' solution branches unique to the LPA-ODE's.  The location and stability of the global branches provides linear stability information for the RDE's \eqref{LPA-1}.  The location and stability of the local branches provides non-linear stability information.  Application of this method and interpretation of its results were demonstrated using two classical example systems, Schnakenberg and substrate inhibition.  To demonstrate scalability of the LPA to larger more complex systems, it was applied to a biologicaly motivated system involving nine interacting chemotaxis regulators.


The Local Perturbation Analysis has a number of benefits.  i) It is simpler to implement than most PDE analysis techniques.  These are notoriously challenging and typically tailored to the particular system being investigated.  ii) With the help of existing software it can be applied to systems with many variables, common in biological applications.  iii)  In this setting, effects of both parameter and structural variations of a reaction network, resulting from cells using the same conserved biochemical machinery in different ways, can be efficiently investigated.   iv) Since the method requires only analysis of ODE's, it should be accessible to a broad base of uses not versed in non-linear PDE techniques.

Beyond these practical considerations, the LPA concisely summarizes a rich set of linear and non-linear information on a single one or two parameter bifurcation diagram.   v) It accurately detects the location of linear instabilities (when diffusivities are sufficiently different).  vi) It detects non-linear patterning regimes which appear to commonly be associated with sub-critical Turing type bifurcations.  vii) In these regimes, it qualitatively characterizes the dependence of response thresholds on reaction parameters.  viii) The global bifurcation structure can be interpreted to provide reasonable hypotheses about the type of pattern that will evolve on longer timescales.  For these reasons, the LPA has the potential to be of use in an array of scientific fields where RDE's arise.

\section*{Acknowledgments}
WRH thanks Leah Edelstein-Keshet and Michael Ward for comments on this manuscript.  This research was partially supported by the NIH grants R01 GM086882 (to Anders E. Carlsson and Leah Edelstein-Keshet) and P50GM76516, and an NSERC discovery grant (to LEK)

\appendix

\section{Schnakenberg asymptotics}\label{Schnak-asym}

This analysis closely follows \cite{Ward-02}.  Consider the Schnakenberg system \eqref{snack-sys} on the interval $[-1,1]$.  In \cite{Ward-02} it was shown that this system exhibits stable spike solutions when $a=0$.  That analysis can be extended to show such spikes in fact exist for all values of $a$ under certain asymptotic conditions.

To begin, define
\begin{equation}
D=\frac{\bar{D}}{\epsilon}, \;\;\;\;\;\; v=\epsilon \bar{v}, \;\;\;\;\;\; u=\frac{\bar{u}}{\epsilon},
\end{equation}
and subsequently drop the $\bar{}$ to yield
\begin{align}\label{snack-sys-app-a}
u_t(x,t)&=a\epsilon-u+u^2v+ \epsilon  u_{xx}\\ \label{snack-sys-app-b}
\epsilon v_t(x,t)&=b-\frac{u^2v}{\epsilon} + D v_{xx} , 
\end{align}
Assuming $D \gg 1/\epsilon$, $v=v_0+\epsilon v_1(x)+...$, and integrating \eqref{snack-sys-app-b}, it can be determined that 
\begin{equation}\label{equ:v0}
2b=\frac{v_0}{\epsilon} \displaystyle \int_{-1}^1 u^2(x,t) \, dx .
\end{equation}
A spike solution of the form $u(x)=u_0 + u_1(x/\epsilon)$ is now sought where $u_0$ and $u_1$ are the outer and inner solutions.  It is assumed $u_0$ is spatially constant.  Collecting terms involving the same powers of $\epsilon$ shows the outer solution is $u_0 \approx a \epsilon$ and the inner solution satisfies
\begin{equation}
u^{''}_{1}(z)-u_1(z)+u_1^2(z) v_0=0
\end{equation}
on $x/\epsilon=z \in [-1, 1]$ with no flux boundary conditions.  The solution to this problem is known (see \cite{Ward-02}) yielding
\begin{equation}
u(x)=u_0+u_1= a \epsilon + \frac{3}{2 v_0} \sech^2 (\frac{x}{2\epsilon}) .
\end{equation}
Integrating the square of this expression and substituting into \eqref{equ:v0} yields $v_0=b/3$.  Unravelling the change of coordinates yields the approximate spike solution for the original problem \eqref{snack-sys} on $[-1,1]$ 
\begin{equation}\label{equ:spike}
u(x) \approx a + \frac{b}{2 \epsilon} \sech^2 (\frac{x}{2\epsilon}), \;\;\;\;\; v(x) \approx \frac{3 \epsilon}{b}.
\end{equation}

So the Schnakenberg system \eqref{snack-sys} in fact produces spike type solutions for all values of $a$ in the limit $\epsilon \to 0$.  This is in agreement with the results of the LPA in Figure \ref{fig:Schnak}a and the progression of the fold bifurcation (where the spike is lost) to $\infty$ as $a \to \infty$ in the bifurcation analysis in Figure \ref{fig:Schnak}c.  Further, the maximum value of $u$ in \eqref{equ:spike} with $a=0$ compares to good precision with the maximum values shown at $a=0$ in Figure \ref{fig:Schnak}c, supporting these results.

\section{Proof of Theorem \ref{Thrm-EigVal}} \label{Thrm-EigVal-proof}

\begin{figure}[htb] 
  \centerline{ 
\includegraphics[scale=.8]{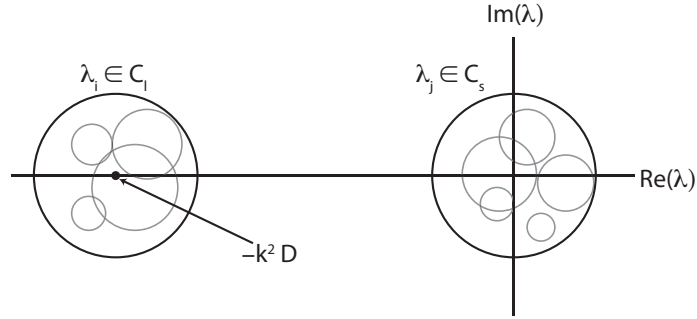}}
\caption{Schematic of the separation of the eigenvalues of $J_k$ in the complex plane.  Grey circles indicate the different Gershgorin circles $C_i$.  The larger darker circles indicate $C_l$ and $C_s$ which contain all eigenvalues of $J_k$.  These circles separate those eigenvalues into two classes with $O(1)$ and $O^-(D)$ real part respectively.  }
\label{fig:Gerhsgorin}
\end{figure} 

To prove Theorem \ref{Thrm-EigVal}, first notice the eigenvalues of $J_k$ \eqref{J_k} can be segregated into two regions of the complex plane using the Gershgorin circle theorem (see Figure \ref{fig:Gerhsgorin}).  Fix a specific wave number $k$, let $a_{i,j}$ be the elements of $J_k$, and define 
\begin{equation}
R_i=\displaystyle \sum_{j \neq i} |a_{i,j}|, \qquad C_i = C(a_{i,i},R_i)
\end{equation}
where $C(a,r)$ is the circle with centre $a$ and radius $r$.  The Gershgorin circle theorem states that each eigenvalue of $J_k$ lies in at least one of the disks $C_i$.  The structure of $J_k$ is such that the off diagonal entries are $O(1)$ with respect to $D$.  So $R=\max \{R_i\}$ is  $O(1)$.  The diagonal entries fall into two categories, those that are $O(1)$ (corresponding to the small diffusion entries), and those that are $-k^2 D + O(1)$ (corresponding to large diffusion entries).  Define $\Omega_s$ to be the union of the disks $C_i$ that are characterized by $O(1)$ diagonal entries and $\Omega_l$ as the union of disks characterized by $O(D)$ diagonal entries.  Since these disks have a maximal radius $R$ independent of $D$, there exists disks $C_l=C(-k^2 D,\kappa_l R)$ and $C_s=C(0, \kappa_s R)$ so that for constants $\kappa_{l,s}$ independent of $D$, $\Omega_s \subset C_s$ and $\Omega_l \subset C_l$.  For sufficiently large $D$, $C_s$ and $C_l$ do not overlap and hence separate $\{\lambda^k\}$ into two sets (see Figure \ref{fig:Gerhsgorin}).  

So, for each $i$, either $\Re(\lambda^k_i)=O^-(D)$, or $\Re(\lambda^k_i)=O(1)$.  Since $\det (J_k) = O(D^{N})$, $\{ \lambda^k_i \}_{i=M+1:M+N}$ must have $O^-(D)$ real part.   Also note that the imaginary parts of all eigenvalues are constrained to be less than $\max \{\kappa_s,\kappa_l \} R$ so that $\Im(\lambda_k^i)=O(1)$ for all $i$ as well, so $|\lambda^k_i|=O(1)$ for $i=1:M$.

Eigenvalues of $J_k$ are roots of the characteristic polynomial 
\begin{equation}
\begin{vmatrix} J_k - \lambda I \end{vmatrix}=\begin{vmatrix} f_u (u^s,v^s) - (k^2 \epsilon^2 + \lambda) I & f_v(u^s,v^s) \\ g_u(u^s,v^s) & g_v(u^s,v^s) - (k^2 D  + \lambda) I \end{vmatrix} = 0,
\end{equation}
where $I$ is a properly sized identity matrix.
Let  $P$ and $Q$ be the unitary matrices that diagonalize $f_u(u^s,v^s)$ and $g_v(u^s,v^s)$.  Then in particular the diagonal entries of $P^{-1} f_u(u^s,v^s) P$ are $\{ \lambda^{LP}_j \}$ and the entries of $Q^{-1} g_v(u^s,v^s) Q$ are $O(1)$.  Define
\begin{equation}
T = \begin{bmatrix} P & 0 \\ 0 & Q \end{bmatrix} .
\end{equation}
Then the eigenvalue problem translates to
\begin{equation}\label{Part-Diag-Jac-EV}
\begin{vmatrix} [\lambda^{LP}] - k^2 \epsilon^2 I - \lambda I & P^{-1} f_v Q \\
 Q^{-1} g_u P & Q^{-1} g_v(u^s,v^s) Q -k^2 D I  - \lambda I 
 \end{vmatrix}
=
\begin{vmatrix}
A1 & A2 \\
A3 & A4
\end{vmatrix}
 =0
\end{equation}
where $[\lambda^{LP}]$ is the diagonal form of $f_u(u^s,v^s)$.  Notice that $A1, A4$ are diagonal.

Now consider an eigenvalue $\lambda$ whose real part is $O(1)$.  In this case, the diagonal entries of $A4$ are $O^-(D)$ and it is non-singular.  It can thus be used to eliminate $A2$.  After this is done, the eigenvalue problem becomes 
\begin{equation}
\begin{vmatrix} [\lambda^{LP}] - k^2 \epsilon^2 I +O(D^{-1})- \lambda I & 0 \\
 Q^{-1} g_u P & Q^{-1} g_v(u^s,v^s) Q - k^2 D I - \lambda I \end{vmatrix}
 =0
\end{equation}
Since the bottom right block is non-singluar, it must be true that 
\begin{equation}\label{perturbed-equ}
\det \left(  [\lambda^{LP}] - k^2 \epsilon^2 I +O(D^{-1})- \lambda I  \right)
 =0 .
\end{equation}
where $O(D^{-1})$ is a properly sized matrix with entries of this size.  With $D=\infty$, the roots of this polynomial are simply $\{ \lambda^{LP}_j - k^2 \epsilon^2 \}$.  It is tempting to view Eqn.~~\eqref{perturbed-equ} as a perturbation of this case and apply some form of perturbation bound.  However, $f_u$ is not Hermitian, which is usually required for such bounds.  Instead, the best we can say is that by continuity of the determinant, the roots of this polynomial satisfy
\begin{equation}
\lambda=\lambda^{LP}_j -k^2 \epsilon^2 + c(D),
\end{equation}
where $c(D) \to 0$ as $D \to \infty$.

\section{GTPase model equations}\label{GTPase-equ}

Figure \ref{fig:GTPase2}a schematically diagrams interactions between three interacting GTPases and three phosphoinositides.  I briefly outline the model equations describing these interactions.  Further specifics can be found in \cite{Holmes12a,Lin-12}.  Modifications of the model presented in those references, which are the subject of investigation here, are described in the main text.  Each GTPase undergoes conservative cycling between active membrane bound and inactive forms in the cell interior by (un)binding to the membrane.  These dynamics are described by
\begin{subequations}\label{equ:GTPase}
\begin{align}\label{equ:crrhonewsupp2}
\frac{\partial G}{\partial t} &= I_G 
\frac{G^c}{G_t} -\delta_G G+D_m G_{xx},  \\ \nonumber
\frac{\partial G^c}{\partial t} &= -I_G 
\frac{G^c}{G_t}+\delta_G G+D_{c} G^c_{xx} ,  
\end{align}
where $G=R,\rho,C$ represents the membrane bound form and $G^c$ represents an inactive cytosolic form .  Phosphoinositides interconvert between three states through the hydrolysis / phosphorylation activity of PI5K, PI3K, PTEN etc. which are not explicitly modeled.  
The GTPase activation rate functions encoding the interactions in Figure \ref{fig:GTPase2}a are defined by
\begin{align}\label{equ:crrhonewsupprhs}
I_C &= \left( \frac{{\hat I}_C}{1+\left(\rho/a_1\right)^n} \right),\quad
I_R = \left( {\hat I}_{R1}+ \frac{\alpha C + \hat{I}_{R2} \frac{P_3}{P_{3b}}}{1+f_2 \left( \rho / a_3 \right)^n }  \right),\\
I_{\rho} &= \frac{{\hat I}_{\rho}}{1+\left(R/a_2\right)^n}. \nonumber
\end{align}
Phosphoinositide kinetics are modeled by linear and mass action kinetics
\begin{align}\label{equ:PIs}
\frac{\partial P_1}{\partial t} &= I_{P1}-\delta_{P1}P_1+k_{21}P_2-f_{PI5K}(R,C,\rho) P_1 + D_P P_{1xx},\\ \nonumber
\frac{\partial P_2}{\partial t} &= -k_{21}P_2+
f_{PI5K}(R,C,\rho) P_1 - 
f_{PI3K}(R,C,\rho) P_2 + 
f_{PTEN}(R,C,\rho) P_3 + D_P P_{2xx},\\ \nonumber
\frac{\partial P_3}{\partial t} &= 
f_{PI3K}(R,C,\rho) P_2 - 
f_{PTEN}(R,C,\rho) P_3 + D_P P_{3xx},
\end{align}
with feedback terms
\begin{align}\label{equ:Rac-PI-feedback}
f_{PI3K}&=\frac{k_{PI3K}}{2} \left(1+\frac{R}{R_t}\right), \quad f_{PI5K}=\frac{k_{PI5K}}{2}\left(1+\frac{R}{R_t}\right) , \\
f_{PTEN}&=\frac{k_{PTEN}}{2} \left(1+\frac{\rho}{\rho_t}\right) . \nonumber
\end{align}
\end{subequations}
See Table \ref{table:params} for a base parameter set for this model.

\begin{table}[!ht]
\begin{center}
\begin{tabular}{|c|l|l|}\hline
Parameter Name & Value &Meaning\\ \hline \hline
$L_0$   & $20 \;\; \mu m$ & Domain size\\ \hline
$C_t, R_t, \rho_t$   & $2.4, 7.5, 3.1 \;\; \mu M$ & Total levels of Cdc42, Rac, and Rho \\ \hline
${\hat I}_c, {\hat I}_{R1}, {\hat I}_{R2}, {\hat I}_{\rho}$  & $2.95, 0.2, 0.2, 6.6 \;\; \mu M s^{-1}$ & Cdc42, Rac, and Rho activation rates\\ \hline
$a_1,a_2, a_3$  & $1.25,1, 1.25 \;\; \mu M$ & Cdc42 and Rho half max inhibition levels\\ \hline
$n$  & $3$ & Hill coefficient for inhibitory connections\\ \hline
$\alpha$ & $0.55 \;\; s^{-1}$ & Cdc42 dependent Rac activation \\ \hline
$\delta_C, \delta_R, \delta_{\rho}$ & $1 \;\; s^{-1}$ & GAP decay rates of activated Rho-proteins\\ \hline

$I_{P1}$ & $10.5 \;\; \mu M/s$ & $\textrm{PIP}_1$ input rate \\ \hline
$\delta_{P1}$ & $0.21 \;\; s^{-1}$ & $\textrm{PIP}_1$ decay rate \\ \hline
$k_{PI5K}, k_{PI3K}, k_{PTEN}$ & $0.084, 0.00072, 0.432 \;\; \mu M^{-1} s^{-1}$ & Baseline conversion rates \\ \hline
$k_{21}$ & $0.021 \;\; s^{-1}$ & Baseline conversion rate \\ \hline
$P_{3b}$ & $0.15 \;\; \mu M $ & Typical level of $\textrm{PIP}_3$\\ \hline

$D_m, D_{c}, D_{P}$ & $0.1, 50, 5 \;\; \mu m^2/s $ & Diffusion Rates \\ \hline
$f_2$ & $1$ & Non dimensional feedback  parameter  \\ \hline

\end{tabular}
\end{center}
\caption{ \textbf{Model Parameters:}
Base parameter set for the model depicted in Figure \ref{fig:GTPase2}a and represented in Eqs.~\eqref{equ:GTPase}.  The primary parameters of interest are $I_{R1}$ which represents a basel activation rate parameter for Rac, and $f_2$ which modulates the strength of inhibitory feedback from Rho to Rac.
}
\label{table:params}
\end{table}

\bibliographystyle{siam}
\bibliography{LPbib}

\begin{thebibliography}{10}

\bibitem{Barrass-06}
{\sc I.~Barrass, E.~J. Crampin, and P.~K. Maini}, {\em Mode transitions in a
  model reaction--diffusion system driven by domain growth and noise}, Bulletin
  of mathematical biology, 68 (2006), pp.~981--995.

\bibitem{Caron-03}
{\sc E.~Caron}, {\em Rac signalling: a radical view}, Nature Cell Biology, 5
  (2003), pp.~185--187.

\bibitem{Crampin-99}
{\sc E.~J. Crampin, E.~A. Gaffney, and P.~K. Maini}, {\em Reaction and
  diffusion on growing domains: scenarios for robust pattern formation},
  Bulletin of mathematical biology, 61 (1999), pp.~1093--1120.

\bibitem{Dawes-07}
{\sc A.~T. Dawes and L.~Edelstein-Keshet}, {\em Phosphoinositides and rho
  proteins spatially regulate actin polymerization to initiate and maintain
  directed movement in a one-dimensional model of a motile cell}, Biophysical
  Journal, 92 (2007), pp.~744--768.

\bibitem{Matcont-03}
{\sc A.~Dhooge, W.~Govaerts, and Yu.~A. Kuznetsov}, {\em Matcont: A matlab
  package for numerical bifurcation analysis of odes}, ACM TOMS, 29 (2003),
  pp.~141--164.

\bibitem{Auto-citation}
{\sc E.~Doedel, A.~Champneys, T.~Fairgrieve, Y.~Kuznetsov, B.~Oldeman,
  R.~Paffenroth, B.~Sandstede, X.~Wang, and C.~Zhang}, {\em Auto-07p:
  Continuation and bifurcation software for ordinary differential equations:
  Continuation and bifurcation software for ordinary differential equations},
  {URL} http://indy. cs. concordia. ca/auto,  (2007).

\bibitem{Doelman-98}
{\sc A.~Doelman, R.A. Gardner, and T.J. Kaper}, {\em Stability analysis of
  singular patterns in the 1d gray-scott model: a matched asymptotics
  approach}, Physica D: Nonlinear Phenomena, 122 (1998), pp.~1 -- 36.

\bibitem{Ferguson-06}
{\sc G.J. Ferguson, L.~Milne, S.~Kulkarni, T.~Sasaki, S.~Walker, S.~Andrews,
  T.~Crabbe, P.~Finan, G.~Jones, S.~Jackson, et~al.}, {\em Pi (3) k$\gamma$ has
  an important context-dependent role in neutrophil chemokinesis}, Nature cell
  biology, 9 (2006), pp.~86--91.

\bibitem{Fu-01}
{\sc Y.~Fu and Z.~Yang}, {\em Rop gtpase: a master switch of cell polarity
  development in plants}, Trends in plant science, 6 (2001), pp.~545--547.

\bibitem{Meinhardt-72}
{\sc A.~Gierer and H.~Meinhardt}, {\em A theory of biological pattern
  formation}, Kybernetik, 12 (1972), pp.~30--39.

\bibitem{Goehring-11}
{\sc N.~W. Goehring, P.~K. Trong, J.~S. Bois, Debanjan Chowdhury, Ernesto~M.
  Nicola, Anthony~A. Hyman, and Stephan~W. Grill}, {\em Polarization of par
  proteins by advective triggering of a pattern-forming system}, Science, 334
  (2011), pp.~1137--1141.

\bibitem{Grieneisen-thesis}
{\sc V.~Grieneisen}, {\em Dynamics of Auxin Patterning in Plant Morphogenesis},
  PhD thesis, University of Utrecht, 2009.

\bibitem{Holmes12b}
{\sc W.~R. Holmes, A.~E. Carlsson, and L.~Edelstein-Keshet}, {\em Regimes of
  wave type patterning driven by refractory actin feedback: Transition from
  static polarization to dynamic wave behaviour}, Phys Biol, 9 (2012),
  p.~046005.

\bibitem{Holmes12a}
{\sc W.~R. Holmes, B.~Lin, A.~Levchenko, and L.~Edelstein-Keshet}, {\em
  Modeling cell polarization driven by synthetic spatially graded rac
  activation}, PLoS Comput Biol, 8 (2012), p.~e1002366.

\bibitem{Huang-03}
{\sc K.C. Huang, Y.~Meir, and N.S. Wingreen}, {\em Dynamic structures in
  escherichia coli: spontaneous formation of mine rings and mind polar zones},
  Proceedings of the National Academy of Sciences, 100 (2003),
  pp.~12724--12728.

\bibitem{Huang-05}
{\sc K.C. Huang and N.S. Wingreen}, {\em Min-protein oscillations in round
  bacteria}, Physical biology, 1 (2005), p.~229.

\bibitem{Iron-00}
{\sc D.~Iron and M.~J. Ward}, {\em A metastable spike solution for a nonlocal
  reaction diffusion model}, SIAM J. Appl. Math, 60 (2000), pp.~778--802.

\bibitem{Iron-04}
{\sc D.~Iron, J.~Wei, and M.~Winter}, {\em Stability analysis of turing
  patterns generated by the schnakenberg model}, Journal of Mathematical
  Biology, 49 (2004), pp.~358--390.

\bibitem{Jilkine-11}
{\sc A.~Jilkine and L.~Edelstein-Keshet}, {\em A comparison of mathematical
  models for polarization of single eukaryotic cells in response to guided
  cues}, PLoS computational biology, 7 (2011), p.~e1001121.

\bibitem{Jilkine-07}
{\sc A.~Jilkine, A.~F.M. Mar\'ee, and L.~Edelstein-Keshet}, {\em Mathematical
  model for spatial segregation of the {Rho-Family GTPases} based on inhibitory
  crosstalk}, Bulletin of Mathematical Biology, 69 (2007), pp.~1943--1978.

\bibitem{Kaper-09}
{\sc H.G. Kaper, S.~Wang, and M.~Yari}, {\em Dynamical transitions of turing
  patterns}, Nonlinearity, 22 (2009), p.~601.

\bibitem{Thomas-79}
{\sc J.~P. Kernevez, G.~Joly, M.~C. Duban, B.~Bunow, and D.~Thomas}, {\em
  Hysteresis, oscillations, and pattern formation in realistic immobilized
  enzyme systems}, Journal of Mathematical Biology, 7 (1979), pp.~41--56.

\bibitem{Wei-05}
{\sc Theodore Kolokolnikov, Michael~J. Ward, and Juncheng Wei}, {\em The
  existence and stability of spike equilibria in the one-dimensional gray scott
  model: The pulse-splitting regime}, Physica D: Nonlinear Phenomena, 202
  (2005), pp.~258 -- 293.

\bibitem{Kolokolnikov-05}
{\sc T.~Kolokolnikov, M.~J. Ward, and J.~Wei}, {\em Pulse-splitting for some
  reaction-diffusion systems in one-space dimension}, Studies in Applied
  Mathematics, 114 (2005), pp.~115--165.

\bibitem{Lewis-93}
{\sc M.A. Lewis and P.~Kareiva}, {\em Allee dynamics and the spread of invading
  organisms}, Theoretical Population Biology, 43 (1993), pp.~141 -- 158.

\bibitem{Li-09}
{\sc F.~Li and W.M. Ni}, {\em On the global existence and finite time blow-up
  of shadow systems}, Journal of Differential Equations, 247 (2009), pp.~1762
  -- 1776.

\bibitem{Lin-12}
{\sc Benjamin Lin, William~R. Holmes, C.~Joanne Wang, Tasuku Ueno, Andrew
  Harwell, Leah Edelstein-Keshet, Takanari Inoue, and Andre Levchenko}, {\em
  Synthetic spatially graded rac activation drives cell polarization and
  movement}, Proceedings of the National Academy of Sciences, Early Edition
  (2012).

\bibitem{Maree-06}
{\sc A.~F.M. Mar\'ee, A.~Jilkine, A.~Dawes, V.~A. Grieneisen, and
  L.~Edelstein-Keshet}, {\em Polarization and movement of keratocytes: A
  multiscale modelling approach}, Bulletin of Mathematical Biology, 68 (2006),
  pp.~1169--1211.

\bibitem{Mori-08}
{\sc Y.~Mori, A.~Jilkine, and L.~Edelstein-Keshet}, {\em Wave-pinning and cell
  polarity from a bistable reaction-diffusion system}, Biophysical Journal, 94
  (2008), pp.~3684--3697.

\bibitem{Mori-11}
\leavevmode\vrule height 2pt depth -1.6pt width 23pt, {\em Asymptotic and
  bifurcation analysis of wave-pinning in a reaction-diffusion model for cell
  polarization.}, SIAM J Applied Math, 71 (2011), pp.~1401--1427.

\bibitem{Murray-82}
{\sc J.D. Murray}, {\em Parameter space for turing instability in reaction
  diffusion mechanisms: A comparison of models}, Journal of Theoretical
  Biology,  (1982), pp.~143--163.

\bibitem{Murray-02}
{\sc J.~D. Murray}, {\em Mathematical biology : An introduction , third
  edition}, Interdisciplinary Applied Mathematics,  (2002).

\bibitem{Nishiura-82}
{\sc Y.~Nishiura}, {\em Global structure of bifurcating solutions of some
  reaction-diffusion systems}, SIAM J. Appl. Math, 13 (1982), pp.~555--593.

\bibitem{Rubinstein-99}
{\sc L.M. Pismen and B.Y. Rubinstein}, {\em Computer tools for bifurcation
  analysis: general approach with application to dynamical and distributed
  systems}, International Journal of Bifurcation and Chaos, 9 (1999),
  pp.~983--1008.

\bibitem{Rodrigues-11}
{\sc L.A.D. Rodrigues, D.C. Mistro, and S.~Petrovskii}, {\em Pattern formation,
  long-term transients, and the turing--hopf bifurcation in a space-and
  time-discrete predator--prey system}, Bulletin of mathematical biology, 73
  (2011), pp.~1812--1840.

\bibitem{Rubinstein-12}
{\sc B.~Rubinstein, B.D. Slaughter, and R.~Li}, {\em Weakly nonlinear analysis
  of symmetry breaking in cell polarity models}, Physical Biology, 9 (2012),
  p.~045006.

\bibitem{Sanders-99}
{\sc E.E. Sander, P.~Jean, S.~Van~Delft, R.A. Van Der~Kammen, and J.G.
  Collard}, {\em Rac downregulates rho activity reciprocal balance between both
  gtpases determines cellular morphology and migratory behavior}, The Journal
  of cell biology, 147 (1999), pp.~1009--1022.

\bibitem{Schnakenberg-79}
{\sc J.~Schnakenberg}, {\em Simple chemical reaction systems with limit cycle
  behaviour}, Journal of Theoretical Biology, 81 (1979), pp.~389--400.

\bibitem{Short-10}
{\sc M.B. Short, A.L. Bertozzi, and P.J. Brantingham}, {\em Nonlinear patterns
  in urban crime: Hotspots, bifurcations, and suppression}, SIAM Journal on
  Applied Dynamical Systems, 9 (2010), pp.~462--483.

\bibitem{Turing-52}
{\sc A.M. Turing}, {\em The chemical basis of morphogenesis}, Philosophical
  Transactions of the Royal Society of London. Series B, Biological Sciences,
  237 (1952), pp.~37--72.

\bibitem{Nishiura-12}
{\sc K.I. Ueda and Y.~Nishiura}, {\em A mathematical mechanism for
  instabilities in stripe formation on growing domains}, Physica D: Nonlinear
  Phenomena, 241 (2012), pp.~37 -- 59.

\bibitem{vanLeeuwen-97}
{\sc F.N. van Leeuwen, H.E.T. Kain, R.A. van~der Kammen, F.~Michiels, O.W.
  Kranenburg, and J.G. Collard}, {\em The guanine nucleotide exchange factor
  tiam1 affects neuronal morphology; opposing roles for the small gtpases rac
  and rho}, The Journal of cell biology, 139 (1997), pp.~797--807.

\bibitem{Ward-02}
{\sc M.~J. Ward and J.~Wei}, {\em The existence and stability of asymmetric
  spike patterns for the schnakenberg model}, Studies in Applied Mathematics,
  109 (2002), pp.~229--264.

\end{thebibliography}

\end{document}